\documentclass[10pt,a4paper]{elsarticle}
\usepackage[utf8]{inputenc}
\usepackage[T1]{fontenc}
\usepackage{amsmath}
\usepackage{amsfonts}
\usepackage{amssymb}
\usepackage{graphicx}
\usepackage{mathtools,amssymb}
\usepackage{subfig}
\usepackage{graphicx}
\usepackage{optidef}
\usepackage{xcolor}
\usepackage{amsthm}
\usepackage{comment}
\usepackage{bm}
\usepackage[ruled]{algorithm2e}
\usepackage{xspace}
\usepackage{ulem}
\usepackage{multirow}
\usepackage{graphicx}
\usepackage{url}

\newtheorem{theorem}{Theorem}[section]

\usepackage[margin=1in]{geometry}

\def\AMD{{\sf AMMDRPG\xspace}}

\definecolor{armygreen}{rgb}{0.19, 0.53, 0.43}
\definecolor{atomictangerine}{rgb}{1.0, 0.6, 0.4}
\newcommand{\JP}[1]{{\color{armygreen}#1}}

\renewcommand{\arraystretch}{1.5}

\title{Coordinating drones with mothership vehicles: The mothership and multiple drones routing problem with Graphs}

\author[1]{Lavinia Amorosi\corref{cor1}}
\ead{lavinia.amorosi@uniroma1.it}
\author[2]{Justo Puerto\corref{cor1}}
\ead{puerto@us.es}
\author[2]{Carlos Valverde\corref{cor1}}
\ead{cvalverde@us.es}

\cortext[cor1]{Equally contributing authors}

\address[1]{Department of Statistical Sciences, Sapienza University of Rome, Italy}
\address[2]{Department of Statistical Sciences and Operational Research, University of Seville, Spain}

\date{\today}

\begin{document}

\begin{abstract}
This paper considers the optimization problems that arise to coordinate a tandem between a mothership vehicle and a fleet of drones. 
Each drone can be launched from the mothership to perform an operation. After completing the operations the drones return to the mothership to recharge batteries and to be ready for a new operation. Operations consist on (partially) visiting graphs with given lengths to deliver some service or to perform some surveillance/inspection activity.  The goal is to minimize the overall weighted distance traveled by the mothership while satisfying some requirements in terms of percentages of visits to the target graphs. In all cases, we develop exact formulations resorting to mixed integer second order cone programs that are compared on a testbed of instances to assess their performance. We also develop a matheuristic algorithm that provides reasonable solutions.  Computational experiments show the usefulness of our methods in different scenarios. 
\end{abstract}


\maketitle

\section{Introduction}
In recent years the grow of the potential business opportunities related to the use of drone technology has motivated the appearance of an interesting body of methodological literature on optimizing of the use of such technology. 
We can find examples of that in many different sectors, like telecommunication where drones can be adopted in place of traditional infrastructures to provide connectivity (see for example \cite{art:Amorosi2018}, \cite{art:Chiaraviglio2018}, \cite{Jimenez2018}, \cite{art:Amorosi2019}, and \cite{art:Chiaraviglio2019a}), or to temporary deal with the damages caused by a disaster (\cite{art:Chiaraviglio2019}), deliveries (see for example \cite{art:Mathew2015} , \cite{art:Ferrandez2016}, \cite{art:Poikonen2020} and \cite{art:Amorosi2020}), also in emergency contexts (\cite{art:Wen2016}), inspection (\cite{art:Trotta2018}) and others.
The reader is referred to the recent surveys \cite{art:Otto2018} and \cite{art:Chung2020} for further details.\\
\noindent
Among the different aspects that can be considered we want to focus, for its relationship to the development in this paper, to the design, coordination and optimization of the combined routes of drones with a base vehicle. After the initial paper \cite{MURRAY201586} by Murray and Chu, where a combined model of truck and drone is considered, the work of Ulmer and Thomas \cite{Ulmer2018} also considers another model where trucks and drones are dispatched as order are placed and analyze the effect of different policies for either the truck or the drone. Other papers, as for instance, \cite{art:Campbell2017}, \cite{art:Carlsson2017} and \cite{art:Dayarian2017}, have also considered hybrid truck-and-drone models in order to mitigate the limited delivery range of drones. Poikonen and Golden, in \cite{Poikonen2019}, advance on the coordination problem considering the \textit{Mothership and drone routing problem} where these two vehicles are used to design a route that visits a number of points allowing the truck to launch and recover the drone in a continuous space. More recently in \cite{art:Poikonen2020} the authors consider the \textit{$k$-Multi-Visit drone routing problem} where a truck that acts as a mobile depot only allowed to stop in a predefined set of points, launches drones that can deliver more than one package to their designated destination points.
\noindent
Many of these papers make the assumptions that the set of allowable locations to launch/retrieve a drone are fixed and known a priori, the operations performed by the drone consist of delivering to a single point and the coordination is between a truck and a single drone. These assumptions may be appropriate in some frameworks but in other cases it may be better to relax them.\\
\noindent
In particular, only few papers in literature focus on drones operations consisting in traversing graphs rather than visiting single points. Campbell et al. 
\cite{art:Campbell2018} introduce the \textit{Drone Rural Postman Problem} (DRPP). The authors present a solution algorithm based on the approximation of curves in the plane by polygonal chains and that iteratively increases the number of points in the polygonal chain where the UAV can enter or leave. Thus, they solve the problem as a discrete optimization problem trying to better define the curve by increasing the number of points. The authors consider also the case in which the drone has limited capacity and thus it cannot serve all the lines. To deal with this latter case, they assume to have a fleet of drones and the problem consists in finding a set of routes, each of limited length.\\
In \cite{art:CAMPBELL202160}  this problem has been defined as the \textit{Length Constrained K-drones Rural Postman Problem}, a continuous optimization problem where a fleet of homogeneous drones have to jointly service (traverse) a set of (curved or straight) lines of a network. The authors design and implement a branch-and-cut algorithm for its solution and a matheuristic algorithm capable of providing good solutions for large scale instances of the problem.\\
Scanning the literature of arc routing problems involving hybrid systems consisting in one vehicle and one or multiple drones, the number of contributions is rather limited.\\
In \cite{art:Tokekar2016} the authors study the path planning problem of a system composed of a ground robot and one drone in precision agriculture and solve it by applying orienteering algorithms. Also the paper \cite{art:Garone2010} studies the problem of paths planning for systems consisting in a carrier vehicle and a carried one to visit a set of target points and assuming that the carrier vehicle moves in the continuous space.\\
To the best of our knowledge, only the paper \cite{art:Amorosi2021}, deals with the coordination of a mothership with one drone to visit targets represented by graphs. The authors made different assumptions on the route followed by the mothership: it can move on a continuous framework (the Euclidean plane), ii) on a connected piecewise linear polygonal chain or iii) on a general graph. In all cases, the authors develop exact formulations resorting to mixed integer second order cone programs and propose a matheuristic algorithm capable to obtain high quality solutions in short computing time.
\\
In this paper we deal with an extension of the problem studied in \cite{art:Amorosi2021} for which we propose a novel truck-and-multi-drones coordination model. We consider a system where a base vehicle (mothership) can stop anywhere in a continuous space and has to support the launch/retrieve of a number of drones that must visit graphs. The contribution on the existing literature is to extend the coordination beyond a single drone to the more cumbersome case of several drones and the operations to traversing graphs rather than visiting single points. In particular, we focus on a synchronous version in which every drone is launched and retrieved in the same stage. We present a mathematical programming formulation, valid inequalities to reinforce it and a matheuristic to deal with large instances of this problem. Moreover, we discuss and show how to extend further the former case to the  asynchronous situation where one assumes that the mothership can retrieve one drone in a different stage from the one in which it has been launched.\\
\noindent
The work is structured as follows: Section 2 provides a detailed description of the problem under consideration and develops a valid mixed integer non linear programming formulations for it. Section 3  provides some valid inequalities that strengthen the formulation and also derives upper and lower bounds on the big-M constants introduced in the proposed formulation. Section 4 presents details of the matheuristic algorithm designed to handle large instances. In Section 5 we report the results obtained testing the formulation and the matheuristic algorithm on different classes of planar graphs in order to assess its effectiveness. Finally, Section 6 concludes the paper.
\section{Problem description and valid formulation}\label{section:desc}

\subsection{Problem description}
In the All terrain Mothership and Multiple Drones Routing Problem with Graphs (\AMD), there is one mothership (the base vehicle) and a fleet of homogeneous drones $\mathcal D$ that have to coordinate among them and with the mothership to perform a number of operations consisting in visiting given percentages of the length of a set of graphs $\mathcal G$. The mothership and the drones travel at constant velocities $v_M$ and $v_D$, respectively. Moreover, it is assumed that each  drone has a limited flying autonomy (endurance) $N^d$, so that once it is launched it must complete the operation and return back to the base vehicle to recharge batteries before the time limit. The base vehicle can move freely on a continuous space and starts at a known location, denoted $orig$ where the mothership and the fleet of drones are ready to depart. Once all the operations are finished the mothership and the drones must return together to a final location, called $dest$. 
\noindent
The set of target graphs $\mathcal G$ to be visited  permits to model real situations like monitoring and inspection activities on portions of networks (roads or wires) where traditional vehicles cannot arrive, due to, for example,  the presence of narrow streets, or because of a natural disaster or a terrorist attack that caused damages on the network. In all these cases, the inspection or monitoring of the drone consists in traversing edges of the network to perform a reconnaissance activity. For this reason we model the targets, to be visited by the drone, as graphs. The operation of visiting a graph can be of two different types: 1) traversing a given percentage of the length of each one of its edges or 2) visiting a percentage of the total length of the network. 


\noindent
In this problem, it is assumed that each graph must be visited by one drone: once the drone assigned to the operation enters the graph, it has to complete the entire operation of traversing this target before to be able to leave the graph to return to the base. Moreover, at each stage the drones must be launched from the mothership at the same point (the launching points have to be determined) and they also must be retrieved at the same point (the rendezvous points also have to be determined). However, this does not mean that the mothership and all drones must arrive at a rendezvous location at the same time: the fastest arriving vehicle may wait for the others at the rendezvous location. Note also that every drone of the fleet does not have to be launched from the current base vehicle location in all the stages because of the capacity constraint. In addition, it is supposed that the cost induced by the drones' trips are negligible as compared to those incurred by the base vehicle. Therefore, the goal is to minimize the overall distance traveled by the mothership. In spite of that, the reader may note that from a theoretical point of view the extension to include in the objective function also the distances traveled by the drones is straightforward and does not increase the complexity of the models and formulations.
\noindent
The goal of the \AMD \ is to find the launching and rendezvous points of the fleet of drones $\mathcal D$ satisfying the visit requirements for the graphs in $\mathcal G$ and minimizing the length of the path traveled by the mothership.\\

\noindent

\subsection{Mixed Integer Non Linear Programming Formulations}\label{Form}
\noindent
In this section we present a MINLP mathematical programming formulation for the \AMD 
that can be used to solve medium size instances of this problem.
\noindent
As mentioned in Section \ref{section:desc}, we assume that the mothership is allowed to move freely in a continuous space that for the sake of presentation we  assume to be $\mathbb R^2$. Here, distances are measured by the Euclidean norm, $\|\cdot\|_2$, although this assumption can be extended to any $l_p$ norm, $1\leq p\leq \infty$ (see \cite{Blanco2017}).
\noindent
In the following, we introduce the parameters or input data that formally describe the problem and that are summarized in Table \ref{table:t1}.

 \begin{table}[!h]
\scriptsize
\centering
\begin{tabular}{ | l | }
\hline
\textbf{Problem Parameters}\\
\hline
$orig$: coordinates of the point defining the origin of the mothership path (or tour).\\
$dest$: coordinates of the point defining the destination of the mothership path (or tour).\\
$\mathcal{G}$: set of the target graphs.\\
$g = (V_g, E_g)$: set of nodes and edges of each target graph $g \in \mathcal{G}$.\\
$\mathcal{L}(e_g)$: length of edge $e$ of graph $g \in \mathcal{G}$.\\
$\mathcal{L}(g)=\sum_{e_g\in E_g} \mathcal L(e_g)$: total length of the graph $g\in\mathcal G$.\\
$B^{e_g}, C^{e_g}$: coordinates of the endpoints of edge $e$ of graph $g \in \mathcal{G}$.\\
$\alpha^{e_g}$: percentage of edge $e$ of graph $g \in \mathcal{G}$ that must be visited.\\
$\alpha^{g}$: percentage of graph $g \in \mathcal{G}$ that must be visited.\\
$v_D$: drone speed.\\
$v_M$: mothership speed.\\
$N^d$: drone endurance. \\
$M$: big-M constant.\\
\hline
\end{tabular}
\caption{Nomenclature for AMMDRPG}
\label{table:t1}
\end{table}

\begin{table}[h!]
\scriptsize
\centering
\begin{tabular}{|l|}
\hline 
\textbf{Binary and Integer Decision Variables}\\
\hline
$\mu^{e_g} \in \{0,1\} \:\: \forall e_g \in E_g$ ($g \in \mathcal{G}$): equal to 1 if edge $e$ of graph $g$ (or a portion of it) is visited by the drone,\\ \hspace*{1cm} and  0 otherwise.\\
$\text{entry}^{e_g} \in \{0,1\} \:\: \forall e_g \in E_g$ ($g \in \mathcal{G}$): auxiliary binary variable used for linearizing expressions.\\
$u^{e_{g}td} \in \{0,1\} \:\: \forall e_g \in E_g$ ($g \in \mathcal{G}$) $\: \forall t \in \mathcal T \: \forall d \in \mathcal D$: equal to 1 if the drone $d$ enters in graph $g$ by the edge $e_g$ at stage $t$,\\ \hspace*{1cm} 0 otherwise.\\
$z^{e_{g}e^{'}_{g}} \in \{0,1\} \:\: \forall e_g, e_g' \in E_g$ ($g \in \mathcal{G}$): equal to 1 if the drone goes from $e_g$ to $e^{'}_{g}$, 0 otherwise.\\
$v^{e_{g}td} \in \{0,1\} \:\: \forall e_g \in E_g$ ($g \in \mathcal{G}$) $\: \forall t \in\mathcal T \: \forall d \in \mathcal D$: equal to 1 if the drone $d$ exits from graph $g$ by $e_g$ at stage $t$,\\ \hspace*{1cm} 0 otherwise.\\
$s^{e_g},\; \forall e_g \in E_g$ ($g \in \mathcal{G}$): integer non negative variable representing the order of visit of the edge $e$ of graph $g$.\\
\hline
\textbf{Continuous Decision Variables}\\
\hline
$\rho^{e_g} \in [0,1]$ and $\lambda^{e_g} \in [0,1] \:\: \forall e_g \in E_g$ ($g \in \mathcal{G}$): defining the entry and exit points on $e_g$.\\
$\nu_\text{min}^{e_g}$ and $\nu_\text{max}^{e_g} \in [0,1] \forall e_g \in E_g$ ($g \in \mathcal{G}$): auxiliary variables used for linearizing expressions.\\
$x_L^t \:\: \forall t \in\mathcal T$: coordinates representing the point where the mothership launches the drone at stage $t$.\\
$x_R^t \:\: \forall t \in\mathcal T$: coordinates representing the point where the mothership retrieves the drone at stage $t$.\\
$R^{e_g} \:\: \forall e_g \in E_g$ ($g \in \mathcal{G}$): coordinates representing the entry point on edge $e$ of graph $g$.\\
$L^{e_g} \:\: \forall e_g \in E_g$ ($g \in \mathcal{G})$: coordinates representing the exit point on edge $e$ of graph $g$.\\
$d_L^{e_gtd} \geq 0, \:\: \forall e_g \in E_g$ ($g \in \mathcal{G}$) $\forall t \in\mathcal T \:\forall d\in\mathcal D$: representing the distance travelled by the drone $d$ from the launching\\
\hspace*{1cm} point $x_L^t$ on the mothership at stage $t$ to the first visiting point $R^{e_g}$ on $e_g$.\\
$d^{e_ge^\prime_g} \geq 0, \:\: \forall e_g, e^\prime_g \in E_g $ ($g \in \mathcal{G}$): representing the distance travelled by the drone from the launching\\
\hspace*{1cm} point $L^{e_g}$ on $e_g$ to the rendezvous point $R^{e^\prime_g}$ on $e^\prime_g$.\\
$d^{e_g} \geq 0, \:\: \forall e_g \in E_g$ ($g \in \mathcal{G}$): representing the distance travelled by the drone from the rendezvous\\
\hspace*{1cm} point $R^{e_g}$ to the launching point $L^{e_g}$ on $e_g$. \\
$d_R^{e_gtd} \geq 0 \:\: \forall e_g \in E_g$ ($g \in \mathcal{G}$) $\forall t \in\mathcal T\:\forall d\in\mathcal D$: representing the distance travelled by the drone $d$ from the last\\
\hspace*{1cm} visiting point $L^{e_g}$ on $e_g$ to the rendezvous point $x_R^t$ on the mothership at stage $t$.\\
$d_{LR}^t \geq 0 \:\: \forall t \in\mathcal T$: representing the distance travelled by the mothership from the \\
\hspace*{1cm} launching point $x_L^t$ to the rendezvous point $x_R^t$ at stage $t$.\\
$d_{RL}^t \geq 0 \:\: \forall t \in\mathcal T$: representing the distance travelled by the mothership from the \\ 
\hspace*{1cm} rendezvous point $x_R^t$ at stage $t$ to the launching point $x_L^{(t+1)}$ at the stage $t+1$.\\
\hline
\end{tabular}
\caption{Decision Variables for AMMDRPG}
\label{table:t2}
\end{table}

\noindent
To represent the movement of the drone within a graph $g\in\mathcal G$, we proceed to introduce some notation related to $g$.
Let $g = (V_g, E_g)$ be a graph in $\mathcal G$ whose total length is denoted by $\mathcal L(g)$. Here, $V_g$ denotes the set of nodes and $E_g$ denotes the set of edges connecting pairs of nodes. Let $e_g$ be the edge $e$ of the graph $g \in G$ and let $\mathcal  L(e_g)$ be its length. Each edge $e_g$ is parameterized by its endpoints $B^{e_g}= (B^{e_g}(x_1), B^{e_g}(x_2))$ and $C^{e_g}= (C^{e_g}(x_1), C^{e_g}(x_2))$ and we can compute its length $\mathcal L(e_g) =\|C^{e_g} -  B^{e_g}\|$.

\noindent
As discussed in Section \ref{section:desc}, we consider two modes of visit to the target graphs $g\in \mathcal{G}$:
\begin{itemize}
    \item Visiting a percentage $\alpha^{e_g}$ of each edge $e_g$ which can be modeled by using the following constraints:
    \begin{equation}\label{eq:alphaE}\tag{$\alpha$-E}
    |\lambda^{e_g} - \rho^{e_g}|\mu^{e_g}\geq \alpha^{e_g}, \quad \forall e_g\in E_g.
    \end{equation}
    \item Visiting a percentage $\alpha^g$ of the total length of the graph:
    \begin{equation}\label{eq:alphaG}\tag{$\alpha$-G}
    \sum_{e_g\in E_g} \mu^{e_g}|\lambda^{e_g} - \rho^{e_g}|\mathcal L(e_g) \geq \alpha^g\mathcal L(g).
    \end{equation}
\end{itemize}

\bigskip
\noindent
In both cases the corresponding constraints are nonlinear. In order to linearize them, we need to introduce a binary variable $\text{entry}^{e_g}$ that determines the traveling direction on the edge $e_g$ as well as the definition of the parameter values $\nu_\text{min}^{e_g}$ and $\nu_\text{max}^{e_g}$ of the access and exit points to that segment. Then, for each edge $e_g$, the absolute value constraint \eqref{eq:alphaE} can be represented by:

\begin{equation}\label{eq:alpha-E}\tag{$\alpha$-E}
 \mu^{e_g}|\rho^{e_g}-\lambda^{e_g}|\geq \alpha^{e_g} \Longleftrightarrow
 \left\{
 \begin{array}{ccl}
  \rho^{e_g} - \lambda^{e_g}                       & =    & \nu_\text{max}^{e_g} - \nu_\text{min}^{e_g},                                     \\
  \nu_\text{max}^{e_g}                         & \leq & 1-{\text{entry}^{e_g}},                                   \\
  \nu_\text{min}^{e_g}                      & \leq & {  \text{entry}^{e_g}},                                        \\
  \mu^{e_g}(\nu_\text{max}^{e_g} + \nu_\text{min}^{e_g} ) & \geq & \alpha^{e_g}.
  \\
 \end{array}
 \right.
\end{equation}

\noindent
The linearization of \eqref{eq:alphaG} is similar to \eqref{eq:alphaE} by changing the last inequality in \eqref{eq:alpha-E} by

\begin{equation}\label{eq:alpha-G}\tag{$\alpha$-G}
\sum_{e_g\in E_g} \mu^{e_g}(\nu_\text{max}^{e_g} + \nu_\text{min}^{e_g})\mathcal L(e_g)\geq \alpha^g\mathcal L(g).
\end{equation}

\noindent
To model this problem, we use stages identified with the order in which the different elements in the problem are visited. Let us denote by $\mathcal T$ the set of stages/tasks that the mothership and the fleet of drones have to carry out. These stages are visits to the different graphs in $\mathcal G$ with the required constraints. A stage $t\in\mathcal T$ is referred to as the operation in which the mothership launches some drones from a taking-off location, denoted by $x_L^t$ and later it takes them back on a rendezvous location $x_R^t$. Here, it is important to realize that both locations $x_L^t$ and $x_R^t$ must be determined in the continuous space where the mothership is assumed to move. Note that $|\mathcal T|\leq|\mathcal G|$, since it is assumed that, for each stage, at least one drone must be launched.
\noindent
For each stage $t\in\mathcal T$, each one of the drones launched from the mothership must follow a path starting from and returning to the mothership, while visiting the required edges of $g$. According to the notation introduced above, we write this generic path in the following form:

$$
x_L^t\rightarrow R^{e_g}\rightarrow L^{e_g}\rightarrow\ldots\rightarrow R^{e^\prime_g}\rightarrow L^{e^\prime_g}\rightarrow \ldots \rightarrow R^{e''_g} \rightarrow x_R^t\rightarrow x_L^{t+1}.
$$

\begin{figure}[h!]
\centering
\includegraphics[width=0.95\linewidth]{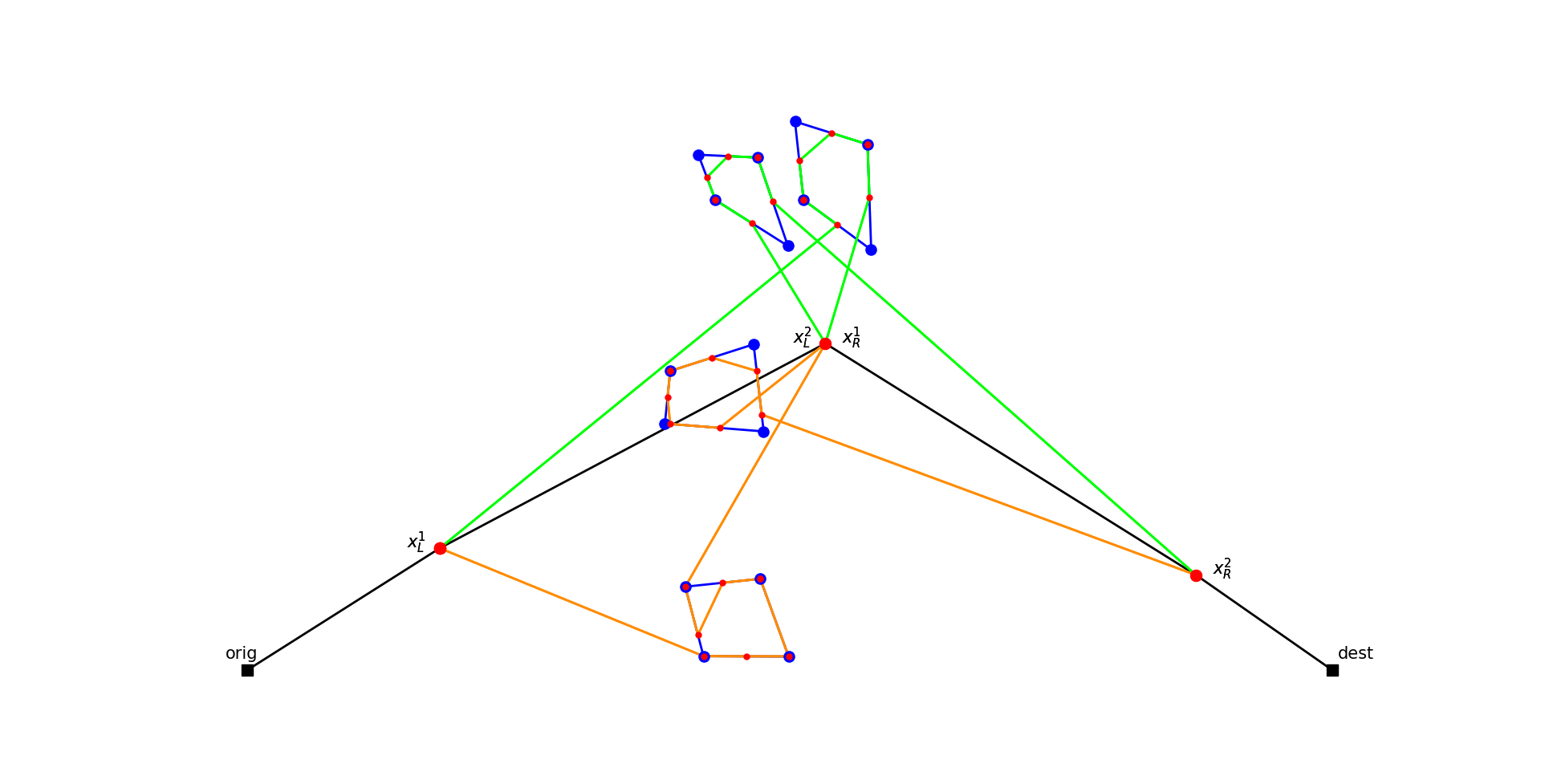}
\caption{Example illustrating the meaning of the launching  (L) and retrieving (R) points.}
\label{fig:illustrative}
\end{figure}

\noindent
Figure \ref{fig:illustrative} shows an example of the notation over a configuration with four target graphs that have four nodes and four edges. Here, it is supposed that the number of available drones is two. 
The mothership begins at its starting point $orig$. Then, it moves to $x_L^1$ where two drones are launched to visit two graphs. There, each drone follows a route (represented by the orange and green paths) that ensures the coverage of $50 \%$ of each edge of the graph. The red dots on the visited graphs are the intermediate points $R^{e_g}$ and $L^{e_g}$ used by the drones in their visit to the edges of the different graphs. After finishing the visit of the first two graphs the drones return to the point $x_R^1$ from where they are launched again to visit the remaining two graphs. Once these graphs are visited, the drones return to the mothership at the rendezvous point $x_R^2$ and then the mothership ends its route at the destination point $dest$.



\noindent
To include the definition of these paths in our mathematical programming formulation we need to make decisions to choose:
\begin{itemize}
    \item The optimal assignment of drones for visiting graphs in a given stage $t$.
    \item The optimal order to visit the edges of each graph in its corresponding stage.
\end{itemize}

\noindent
We model the route that the drone follows by using the binary variables $u^{e_gtd}$, $z^{e_ge^\prime_g}$ and $v^{e_gtd}$ defined in Table \ref{table:t2}.


\begin{align}
    \sum_{g\in \mathcal G}\sum_{e_g\in E_g} \sum_{d\in\mathcal D} u^{e_gtd} & \leq 1, &\forall t\in \mathcal T, \label{st:DEnt}\\
    \sum_{g\in\mathcal G}\sum_{e_g\in E_g} \sum_{d\in\mathcal D} v^{e_gtd} & \leq 1, &\forall t\in \mathcal T, \label{st:DExt}\\
    \sum_{e_g\in E_g} \sum_{t\in \mathcal T} \sum_{d\in\mathcal D} u^{e_gtd} & = 1, &\forall g\in\mathcal G, \label{st:DEng}\\
    \sum_{e_g\in E_g} \sum_{t\in \mathcal T} \sum_{d\in\mathcal D} v^{e_gtd} & = 1, &\forall g\in\mathcal G, \label{st:DExg}\\
    \sum_{e_g\in E_g} u^{e_gtd} & = \sum_{e_g\in E_g} v^{e_gtd}, &\forall g\in\mathcal G, \forall t\in \mathcal T, \forall d\in\mathcal D, \label{st:Duv}\\
     \sum_{t\in \mathcal T} \sum_{d \in \mathcal D} u^{e_gtd} + \sum_{e^\prime_g\in E_g} z_g^{e^\prime_ge_g} & = \mu^{e_g}, &\forall e_g\in E_g:g\in\mathcal G, \label{st:DInu}\\
     \sum_{t\in \mathcal T} \sum_{d \in \mathcal D} v^{e_gtd} + \sum_{e^\prime_g\in E_g} z_g^{e_ge^\prime_g} & = \mu^{e_g}, &\forall e_g\in E_g:g\in\mathcal G. \label{st:DInv}
\end{align}

\noindent 
Inequalities \eqref{st:DEnt} and \eqref{st:DExt} state that for each stage at most one drone can be launched and retrieved for performing an operation.  Constraints \eqref{st:DEng} and \eqref{st:DExg} assure that each graph is visited at some stage $t$ by some drone $d$. Equations \eqref{st:Duv} ensure that the operation of entering and exiting from the graph $g$ occurs in the same stage $t$ and is done by the same drone $d$. Constraints \eqref{st:DInu} state that 
if edge $e$ of graph $g$ is visited by the drone $d$, one of two alternative situations must occur: either $e$ is the first edge of graph $g$ visited by the drone $d$ at stage $t$, or edge $e$ is visited by the drone $d$ after visiting another edge $e^\prime$ of graph $g$. Similarly, constraints \eqref{st:DInv} state that if edge $e$ of graph $g$ is visited by the drone $d$, either $e$ is the last edge of graph $g$ visited by the drone at stage $t$, or the drone $d$ must move to another edge $e^\prime$ of graph $g$ after visiting edge $e$.

\subsubsection*{Elimination of subtours}
\noindent
In order to represent actual routes of the drones over the target graphs, subtours cannot be allowed. Note that subtours would represent fake operations since they would allow free jumps of the drone between different routes at no time. To prevent the existence of subtours within each graph $g\in \mathcal G$ that the drone must visit, one can include, among others, either the compact formulation that uses the Miller-Tucker-Zemlin constraints (MTZ) or the subtour elimination  constraints (SEC).\\
\noindent
For the MTZ formulation, we use the continuous variables $s^{e_g}$, defined in Table \ref{table:t2}, that state the order to visit the edge $e_g$ and set the following constraints for each $g\in\mathcal G$:

\begin{align}
    s^{e_g} - s^{e^\prime_g} + |E_g|z^{e_ge^\prime_g} & \leq |E_g| - 1  , &\quad\forall e_g \neq e_g'\in E_g, \tag{MTZ$_1$} \label{MTZ1}\\
    0 & \leq s^{e_g} \leq |E_g| - 1, &\quad\forall e_g\in E_g.\tag{MTZ$_2$}\label{MTZ2}
\end{align}

\noindent
Alternatively, we can also use the family of subtour elimination constraints for each $g\in\mathcal G$:
\begin{equation}\tag{SEC}\label{SEC}
    \sum_{e_g, e^\prime_g \in S} z_g^{e_ge^\prime_g} \leq |S| - 1, \quad \forall S\subset E_g.
\end{equation}

\noindent
Since there is an exponential number of SEC constraints, when we implement this formulation we need to perform a row generation procedure including constraints whenever they are required by a separation oracle. To find SEC inequalities, as usual, we search for disconnected components in the current solution. Among them, we choose the shortest subtour found in the solution to be added as a lazy constraint to the model.\\

\noindent
 The goal of the \AMD\xspace is to find a feasible solution that minimizes the total distance traveled by the mothership. To account for the different distances among the decision variables of the model we need to set the continuous variables $d_L^{e_gtd}$, $d^{e_g}$, $d^{e_ge^\prime_g}$, $d_R^{e_gtd}$, $d_{RL}^t$ and $d_{LR}^t$, defined in Table \ref{table:t2}. This can be done by means of the following constraints:
 

\begin{align*}
\|x_L^t- R^{e_g}\| & \leq  d_L^{e_gtd},  &\quad \forall e_g:g\in \mathcal{G}, \forall t\in \mathcal T, \forall d\in\mathcal D, \tag{DIST$_{1}$-t} \label{eq:d1}\\
\|R^{e_g}- L^{e_g}\| & \leq  d^{e_g},  &\quad \forall e_g:g\in \mathcal{G}, \tag{DIST$_{2}$-t} \label{eq:d2}\\
\|R^{e_g}- L^{e^\prime_g}\| & \leq  d^{e_ge^\prime_g}, &\quad \forall e_g\neq e_g'\in E_g:g\in \mathcal{G}, \tag{DIST$_{3}$-t} \label{eq:d3}\\
\|L^{e_g}- x_R^t\| & \leq  d_R^{e_gtd}, &\quad \forall e_g:g\in \mathcal{G},\forall t\in T, \forall d\in\mathcal D, \tag{DIST$_{4}$-t} \label{eq:d4}\\
\|x_R^t- x_L^{t+1}\| & \leq  d_{RL}^t, & \quad \forall t\in \mathcal T, \tag{DIST$_{5}$-t} \label{eq:d5}\\
\|x_L^t- x_R^t\| & \leq  d_{LR}^t, & \quad \forall t\in \mathcal T. \tag{DIST$_{6}$-t} \label{eq:d6}\\
\end{align*}

\noindent
The coordination between the drones and the mothership must ensure that the time spent by the drone $d$ to visit the graph $g$ at the stage $t$ is less than or equal to the time that the mothership needs to move from the launching point to the retrieving point during the stage $t$. To this end, we need to define the following coordination constraint for each graph $g\in \mathcal G$, stage $t\in \mathcal T$ and drone $d\in\mathcal D$:

\begin{equation}\tag{DCW}\label{DCW}
\frac{1}{v_D}\left(\sum_{e_g\in E_g} u^{e_gtd}d_L^{e_gtd} + \sum_{e_g, e^\prime_g\in E_g}z^{e_ge^\prime_g}d^{e_ge^\prime_g} + \sum_{e_g\in E_g} \mu^{e_g}d^{e_g} + \sum_{e_g\in E_g} v^{e_gtd}d_R^{e_gtd}\right) \leq \frac{d_{LR}^t}{v_M} + M(1 - \sum_{e_g\in E_g} u^{e_gtd}).
\end{equation}

\noindent
Eventually, we have to impose that the tour of the mothership, together with the drones, starts from the origin $orig$ and ends at the destination $dest$. To this end, we define the following constraints:

\begin{align*}
x_L^0 & =  orig,  \tag{ORIG$_1$} \label{eq:O1} \\
x_R^0 & =  orig,  \tag{ORIG$_2$} \label{eq:O2} \\
x_L^{|\mathcal{G}|+1} & =  dest,  \tag{DEST$_1$} \label{eq:D1} \\
x_R^{|\mathcal{G}|+1} & =  dest.  \tag{DEST$_2$} \label{eq:D2} 
\end{align*}

\noindent
Note that, since the objective function of this problem minimizes the right-hand-side of \eqref{DCW}, this constraint will become an equality and we can model the time capacity constraint for a particular stage $t\in \mathcal T$ by limiting the distance traveled by the mothership for this task $t$:

\begin{equation}\tag{Capacity}\label{CAP}
    d_{LR}^t \leq N^d.
\end{equation}



\noindent
Therefore, putting together all the constraints introduced before, the following formulation minimizes the overall distance traveled by the mothership ensuring the coordination with the fleet of drones while guaranteeing the required coverage of the target graphs.
\begin{mini*}|s|
 {}{\sum_{t\in \mathcal T} (d_{RL}^t + d_{LR}^t)}{}{} \label{AMMDRPG} \tag{AMMDRPG}
 \addConstraint{\eqref{st:DEnt}-\eqref{st:DInv}}{}{}
 \addConstraint{\eqref{MTZ1}-\eqref{MTZ2}} \text{ or }  \eqref{SEC}
 \addConstraint{\eqref{eq:alpha-E} \text{ or } \eqref{eq:alpha-G}}{}{}
 \addConstraint{\eqref{DCW}}{}{}
 \addConstraint{\eqref{CAP}}{}{}
 \addConstraint{\eqref{eq:d1}-\eqref{eq:d6}}{}{}
 \addConstraint{\eqref{eq:O1}-\eqref{eq:D2}.}{}{}
\end{mini*}

\noindent
The objective function accounts for the distances traveled by the mothership. Constraints \eqref{st:DEnt}-\eqref{st:DInv} models the route followed by the drone $d\in\mathcal D$, \eqref{MTZ1} - \eqref{MTZ2} \text{ or } \eqref{SEC} ensure that the displacement of the drone $d\in\mathcal D$ assigned to the target graph $g\in\mathcal G$ is a route, \eqref{eq:alpha-E} \text{ or } \eqref{eq:alpha-G} defines what is required in each visit to a target graph. Finally, constraints (\ref{eq:d1})-(\ref{eq:d6}) set the variables $d_L^{e_gtd}$, $d^{e_g}$, $d^{e_ge^\prime_g}$, $d_R^{e_gtd}$, $d_{RL}^t$ and $d_{LR}^t$, defined in Table \ref{table:t2}, which represent Euclidean distances needed in the model. \\

\noindent
Note that, to deal with the bilinear terms of \eqref{DCW}, we use McCormick's envelopes to linearize them by adding variables $p\geq 0$  representing the products and introducing the following constraints:
\begin{align*}
    p & \leq  M z, \\
    p & \leq  d, \\
    p & \geq m z, \\
    p & \geq d - M(1 - z),
\end{align*}
where $m$ and $M$ are, respectively, the lower and upper bounds of the distance variable $d$. These bounds will be adjusted for each bilinear term in Section \ref{bounds}.

\subsection{The \AMD\xspace without synchronisation}\label{amdasyn}
\noindent
In the \eqref{AMMDRPG} formulation, we assume that every drone is launched and retrieved in the same stage. In this subsection, we show how this assumption can be relaxed.
We consider a variant of the model presented in Section \ref{Form}, in which we assume that the mothership can retrieve one drone in a stage different from the one in which it has been launched. That is, the mothership can move to another point to launch a new drone without having  retrieved the one that was launched before.



\noindent
To deal with this extension, we do not need to define new variables since it is possible to use the same variables that were used in the previous model. First of all, constraint \eqref{st:Duv} must be changed to:
\begin{equation}\label{constraint:Duv-S}
    \sum_{e_g\in E_g} u^{e_gtd} -  \sum_{e_g\in E_g} \sum_{t'\geq t} v^{e_gt'd}=0, \quad\forall  g\in\mathcal G,\forall t\in\mathcal T, \forall d\in\mathcal D.
\end{equation}
This equation states that, if the graph $g$ is assigned to the drone $d$ at the stage $t$, there is another stage $t' \geq t$ in which the drone must come back to the mothership.
\noindent
In addition, if the drone $d$ is launched at $t_1$ and retrieved at $t_2$, this drone cannot be used in any intermediate stage. Hence, for each graph $g\in\mathcal G$ and drone $d\in\mathcal D$ the following constraints must be satisfied:
\begin{equation}\label{constraint:u-S}
 \sum_{e_g\in E_g}\sum_{t=t_1+1}^{t_2} u^{e_gtd}\leq M(2-\sum_{e_g\in E_g} u^{e_gt_1d} - \sum_{e_g\in E_g}v^{e_gt_2d}),\quad t_1 < t_2,
\end{equation}
\begin{equation}\label{constraint:v-S}
 \sum_{e_g\in E_g}\sum_{t=t_1}^{t_2-1} v^{e_gtd}\leq M(2-\sum_{e_g\in E_g} u^{e_gt_1d} - \sum_{e_g\in E_g}v^{e_gt_2d}),\quad t_1< t_2.
\end{equation}
\noindent
Moreover, the coordination constraint \eqref{DCW} must be modified to consider the general case in which the stages of launching and retrieving can be different. For each $t_1<t_2$ and $\forall g\in\mathcal G,\forall d\in\mathcal D$:
\begin{tiny}
\begin{align}\tag{DCW-S}\label{constraint:DCW-S}
\frac{1}{v_D}\left(\sum_{e_g\in E_g} u^{e_gt_1d}d_L^{e_gt_1d} + \sum_{e_g, e^\prime_g\in E_g}z^{e_ge^\prime_g}d^{e_ge^\prime_g} + \sum_{e_g\in E_g} \mu^{e_g}d^{e_g} + \sum_{e_g\in E_g} v^{e_gt_2d}d_R^{e_gt_2d}\right) \leq & \frac{\sum_{t=t_1}^{t_2}d_{LR}^t}{v_M} + \frac{\sum_{t=t_1}^{t_2-1}d_{RL}^t}{v_M} \nonumber
\\  &  + M(2 - \sum_{e_g\in E_g} u^{e_gt_1d} - \sum_{e_g\in E_g} v^{e_gt_2d}).  \nonumber
\end{align}
\end{tiny}

\begin{figure}[h!]
    \centering
    \includegraphics[width = 0.7\linewidth]{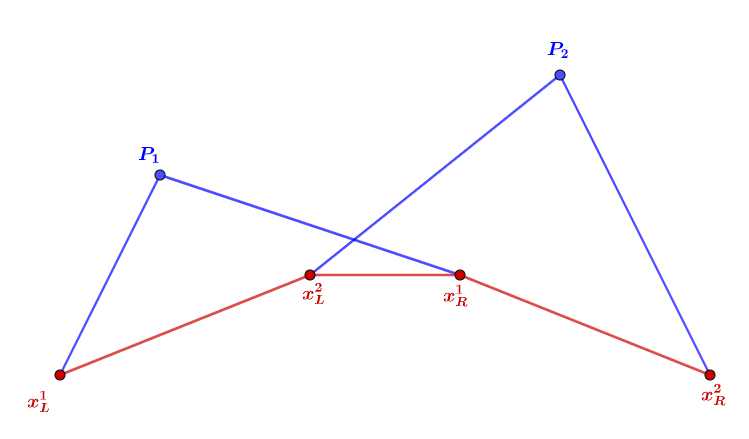}
    \caption{The mothership launches two drones sequentially}
    \label{fig:proof1}
\end{figure}

\noindent
This inequality takes into account the total time the mothership needs to go from the launching point $x_L^{t_1}$ to the rendezvous point $x_R^{t_2}$ when  the drone $d$ is chosen to visit the graph $g\in\mathcal G$ and this operation begins at stage $t_1$ and ends at $t_2$.
\noindent
Finally, we present a result that links the two models presented before. Note that the only difference that solutions can have between these models is that, for the non-synchronized case, the mothership can launch a second drone sequentially before retrieving another one that was launched before. Figure \ref{fig:proof1} shows a solution that is not possible for the model with synchronization. Indeed, we can see that a first drone is launched at $x_L^1$ to visit $P_1$ that is retrieved at $x_R^1$. However, the mothership has launched another drone at $x_L^2$ that goes visiting $P_2$ before having retrieved the first drone. Obviously, this solution does not satisfy the assumption in the synchronized model.

\noindent
In order to present our next result, wlog,  we restrict ourselves to the degenerate case where graphs reduce to points. The reader may note that we could reduce the general case to this one reducing the available capacity so that is possible to traverse the required percentage of these graphs. We simplify the proof considering a generic solution between two consecutive target points.

\begin{theorem}
Let $x_L^1$, $x_L^2$ (resp. $x_R^1$, $x_R^2$) be the launching (resp. rendezvous) points associated to the visit of the target points $P_1$ and $P_2$. If there exist two points $x_L$ and $x_R$ verifying 
$$
 \left\{
 \begin{array}{ccl}
  \dfrac{\|x_L-x_R\|}{v_C} & \leq    & \dfrac{\|x_L - P_1\| + \|P_1 - x_R\|}{v_D}, \\
  \dfrac{\|x_L-x_R\|}{v_C} & \leq    & \dfrac{\|x_L - P_2\| + \|P_2 - x_R\|}{v_D}, \\
  \dfrac{\|x_L-x_R\|}{v_C} & \leq   & N^d, \\
  \|x_L-x_R\| & \leq & \|x_L^1 - x_L^2\| + \|x_L^2- x_R^1\| + \|x_R^1-x_R^2\|,
 \end{array}
 \right.
$$

\noindent then the contribution of this partial route to the optimal objective value will be the same in both models.
\end{theorem}

\begin{proof}
Note that in the considered configuration, the order of visit to the points $P_1$ and $P_2$ is fixed and then, the binary variables in the model are fixed in this case. Thus, the only difference that the two models can have are the location of the launching and rendezvous points. Hence, the only constraints that are involved are those related to these points. These are the conditions in the statement: The first two are the \eqref{DCW} inequalities. The third one is the \eqref{CAP} constraint and the last one ensures that the distance traveled by the mothership in the synchronized model is smaller than or equal to the distance assumed in the non-synchronized solution described in the statement. Therefore, the conclusion follows.

\end{proof}




\section{Strengthening the formulation of \AMD}\label{bounds}
\noindent
In this section we present some valid inequalities for \AMD\xspace that reinforce the formulation given in Section \ref{Form}. Moreover, the \eqref{DCW} constraint has products of binary and continuous variables that, when they are linearized, produce big-M constants that have to be tightened. This section also provides some bounds for these constants whenever it is possible. 



\subsection{Valid inequalities for the \AMD}
\noindent
In this problem, we assume that the fleet has more than one drone since otherwise the problem reduces to (AMDRPG) that was already studied in \cite{art:Amorosi2021}. Therefore, if there exists a stage in which more than one drone is launched, the mothership does not need to perform $|\mathcal G|$ different stages. Hence, most likely the model does not need to deal with those stages that are numbered at the end. By exploiting this idea, it is possible to concentrate all drone operations on the first stages, avoiding empty tasks in $\mathcal T$.
\noindent
Let $\beta^t$ be a binary variable that assumes the value one if all the target graphs are visited when the operation $t$ begins, and  zero, otherwise. Note that, if all the operations are completed before the stage $t$ then they are also completed before the stage $t+1$. Hence, $\beta$ variables must satisfy the following constraints:

\begin{equation}\tag{Monotonicity}\label{eq:Monotonicity}
\beta^t \leq \beta^{t+1}, \mbox{ for all } t=1,\ldots, |\mathcal{G}|-1.
\end{equation}

\noindent
Let $k^t$ denote the number of graphs that are visited in the stage $t$. This number can be computed using the $u$ variables since $u^{e_gtd}$ takes the value 1 if the graph $g$ is visited in stage $t$ by the drone $d$. Thus:

$$k^t=\sum_{e_g\in g:g\in\mathcal G}\sum_{d\in\mathcal D} u^{e_gtd}.$$

\noindent
Hence, if $\beta^t$ equals one, the entire set of graphs in $\mathcal G$ must have been visited before the stage $t$:

\begin{equation}\tag{VI-1}\label{eq:VI-1}
\sum_{t'=1}^{t-1} k^{t'} \geq |\mathcal G|\beta^t,
\end{equation}
where $|\mathcal G|$ denotes the number of graphs of $\mathcal G$.

\noindent
To reduce the space of feasible solutions, we can assume without loss of generality that it is not permitted to have a stage $t$ without any operation if some graphs are still to be visited. This can be enforced by the following constraints:

\begin{equation}\tag{VI-2}\label{eq:VI-2}
k^t \geq 1 - \beta^t.
\end{equation}

\noindent
In addition, it is also possible to reduce the symmetry. Since we are assuming that drones are indistinguishable, we can assume that given an arbitrary order on them, we always assign drones to operations in that given order. This assumption allows us to assign for an operation at the stage $t$ the first drone that is available, avoiding to consider the last ones, if they are not necessary. This consideration can be implemented by means of the following set of inequalities. For all $t\in\mathcal T$:
\medskip

\begin{equation}\tag{VI-3}\label{eq:VI-3}
\sum_{e_g\in \mathcal G} u^{e_gtd} \leq \sum_{e_g:g\in\mathcal G}u^{e_gtd-1}, \; \forall d=2,\ldots |\mathcal D|,      
\end{equation}
\begin{equation}\tag{VI-4}\label{eq:VI-4}
\sum_{e_g\in \mathcal G} v^{e_gtd} \leq \sum_{e_g:g\in\mathcal G}v^{e_gtd-1}, \; \forall d=2,\ldots |\mathcal D|.      
\end{equation}

\noindent
Hence, if the drone $d_1$ is assigned to the task $t$, every drone $d_2$ that is, for the launching order in $\mathcal D$, before than $d_1$, must have been also assigned to perform the task $t$.
\medskip

\noindent
The model that we have proposed includes big-M constants. We have defined different big-M constants along this work. In order to strengthen the formulations we provide tight upper and lower bounds for those constants. In this section we present some results that adjust them for each one of the models.

\subsubsection*{Big $M$ constants bounding the distance from the launching / rendezvous point on the path followed by the mothership to the rendezvous / launching point on the target graph $g\in \mathcal{G}$}

\noindent
To linearize the first term of the objective function in \AMD, we define the auxiliar non-negative continuous variables $p_L^{e_gtd}$ (resp. $p_R^{e_gtd}$) and we model the product by including the following constraints:
\begin{align*}
p_L^{e_gtd} & \geq m_L^{e_g} u^{e_gtd}, \\
p_L^{e_gtd} & \leq d_L^{e_g} - M_L^{e_gtd}(1-u^{e_gt}).
\end{align*}
The best upper bound $M_L^{e_gtd}$ or $M_R^{e_gtd}$ that we can consider is the full diameter of the data, that is the maximum distance between every pair of vertices of the graphs $g\in \mathcal{G}$, in the input data, i.e., every launching or rendezvous point is inside the circle whose diametrically opposite points are described below. 
$$
M_R^{e_gtd} = \max_{\{v\in V_g, v'\in V_{g'} : g, g'\in\mathcal G\}} \|v - v'\| = M_L^{e_gtd}.
$$
\noindent
On the other hand, the minimum distance in this case can be zero. This bound is attainable whenever the launching or the rendezvous points of the mothership is the same that the rendezvous or launching point on the target graph $g\in \mathcal{G}$.




\subsubsection*{Bounds on the big $M$ constants for the distance from the launching to the rendezvous points on the target graph $g\in \mathcal{G}$.} 
\noindent
When the drone visits a graph $g$, it has to go from one edge $e_g$ to another edge $e'_g$ depending on the order given by $z^{e_ge_g'}$. This fact produces a product of variables linearized by the following constraints:
\begin{align*}
p^{e_ge'_g} & \geq m^{e_ge_g'} d_{RL}^{gg'}, \\
p^{e_ge_g'} & \leq d^{e_ge_g'} - M^{e_ge_g'}(1-z^{e_ge_g'}).
\end{align*}

\noindent
Since we are taking into account the distance between two edges $e=(B^{e_g},C^{e_g}), \, e'=(B^{e^\prime_g},C^{e^\prime_g})\in E_g$, the maximum and minimum distances between their vertices give us the upper and lower bounds:
\begin{align*}
M^{e_g e^\prime_g} = & \max\{\|B^{e_g} - C^{e^\prime_g}\|, \|B^{e_g} - B^{e^\prime_g}\|, \|C^{e_g} - B^{e^\prime_g}\|, \|C^{e_g} - C^{j_g}\|\}, \\
m^{e_g e^\prime_g} = & \min\{\|B^{e_g} - C^{e^\prime_g}\|, \|B^{e_g} - B^{e^\prime_g}\|, \|C^{e_g} - B^{e^\prime_g}\|, \|C^{e_g} - C^{e^\prime_g}\|\}.
\end{align*}


\subsubsection*{Bounds on the big $M$ constants for the distance covered by the drone during an operation for all the models by stages}
\noindent
To link the drone operation with the trip followed by the mothership, we have defined the constraint $\eqref{DCW}$ that includes another big-M constant:
\begin{equation*}
\left(\sum_{e_g\in E_g} u^{e_gtd}d_L^{e_gtd} + \sum_{e_g, e^\prime_g\in E_g}z^{e_ge^\prime_g}d^{e_ge^\prime_g} + \sum_{e_g\in E_g} \mu^{e_g}d^{e_g} + \sum_{e_g\in E_g} v^{e_gtd}d_R^{e_gtd}\right)/v_D \leq d_{RL}^t/v_M + M(1 - \sum_{e_g\in E_g} u^{e_gtd}).
\end{equation*}
\noindent
To obtain an upper bound on $M$ we add to the length of the graph $\mathcal L(g)$ the big-Ms computed for $u^{e_gtd}$ and $v^{e_gtd}$, namely $M_{L}^{e_gtd}$ and $M_R^{e_gtd}$, respectively, and the maximum distance that can be traveled by the drone to move from one edge to another one. This results in a valid value for this $M$ constant:

$$M = \mathcal{L}(g) + M_L^{e_gtd} + M_R^{e_gtd} + \sum_{e_g, e_g'\in E_g}M^{e_ge_g'}.$$

\section{A Matheuristic for the Mothership-Drone Routing Problem with Graphs}\label{Math}
\noindent
This section is devoted to present our matheuristic approach to address the solution of the \AMD. Our motivation comes from the fact that the exact solution of the models presented in Section \ref{Form} is highly time demanding. Alternatively, the matheuristic provides good quality solution in limited computing times.\\
\noindent
The basic idea of the algorithm is to determine the route that a drone should perform for visiting each graph $g \in \mathcal{G}$, and thus the entry and exit points $L^{e_{g}}$ and $R^{e^{'}_{g}}$ for each graph.
Sequentially, a clustering procedure on the target graphs is applied in order to compute the route of the mothership via their reference points and the origin/destination points.
The clustering procedure is based on a random selection of the initial target graphs and for this reason it is repeated a number of times in order to consider different cluster structures. At each iteration the new clusters are evaluated by computing the cost of the route visiting their reference points and the origin/destination points. 
The route computed on the reference points of the best cluster generated by this iterative procedure, is used to set the values of the binary variables $u^{e_gtd}$ and $v^{e_gtd}$, that determine the order of visits to the graphs. Finally, these variables are provided as an initial partial solution to the \AMD\xspace model to produce a complete feasible solution.\\
In the following, we present the pseudo-code of this algorithm:

\begin{itemize} 
\item[STEP 1] (First entry and last exit points for each target graph)\\
Compute the route on each target graph $g \in \mathcal{G}$.
Let $L^{e_{g}}$ and $R^{e^{'}_{g}}$ be the pair of entry and exit points on $g$ closest to the origin and let  $\mathcal L(e_{g}, e^{'}_{g})$ be the associated length computed as the sum of the distances travelled by the drone to visit the graph $g$, excluding the distance between $L^{e_{g}}$ and $R^{e^{'}_{g}}$.
\item[STEP 2] (Clustering procedure)\\
Initialization: set $it=1$, define one cluster for each target graph and set $nit=1$. \\
Select randomly two clusters $K_i$ and $K_j$ (where $i<j$).\\
Check if the number of graphs belonging to the union of $K_i$ and $K_j$ is less than the number of available drones $n_D$.\\
If this condition is satisfied:\\
search for point $P$ satisfying the following capacity constraint:
$$
\frac{d(P, R^{e_g}) + \mathcal L(e_{g}, e^{'}_{g}) + d(L^{e^{'}_{g}}, P)}{v_D} \leq N^d, \quad \forall R^{e_g}, L^{e^{'}_{g}} \in K_i, \:\: K_j.
$$
If such a point exists, merge the two clusters and label the new one as $K_i$.\\
Set $nit=nit+1$.\\
Repeat the same procedure on the new cluster structure while $nit < maxit$.
\item[STEP 3] (Computation of Reference Points) 
Compute a reference point for each cluster generated at STEP 2. This computation seeks for the minimization of the distance between each pair of reference points and the distance between them and the origin, always imposing that the \eqref{CAP} constraint is satisfied.

\item [STEP 4] (Setting the order of visits to the  graphs: route via the reference points and the origin/destination points) \\
Compute the TSP of the mothership among the reference points of the clusters and let $\mathcal L(TSP)$ be the associated length.\\
Set $it=it+1$.\\
if(it< maxseed) go to STEP 2\\
else go to STEP 5
\item [STEP 5] (Solution of the \AMD\space model fixing an initial partial solution)\\
Set the values of the binary variables $u^{e_{g}td}$ and $v^{e_{g}td}$ and solve the model \AMD\, to obtain a feasible solution.
\end{itemize}
\noindent

\begin{figure}[h]
\centering
\includegraphics[width = 0.5\linewidth]{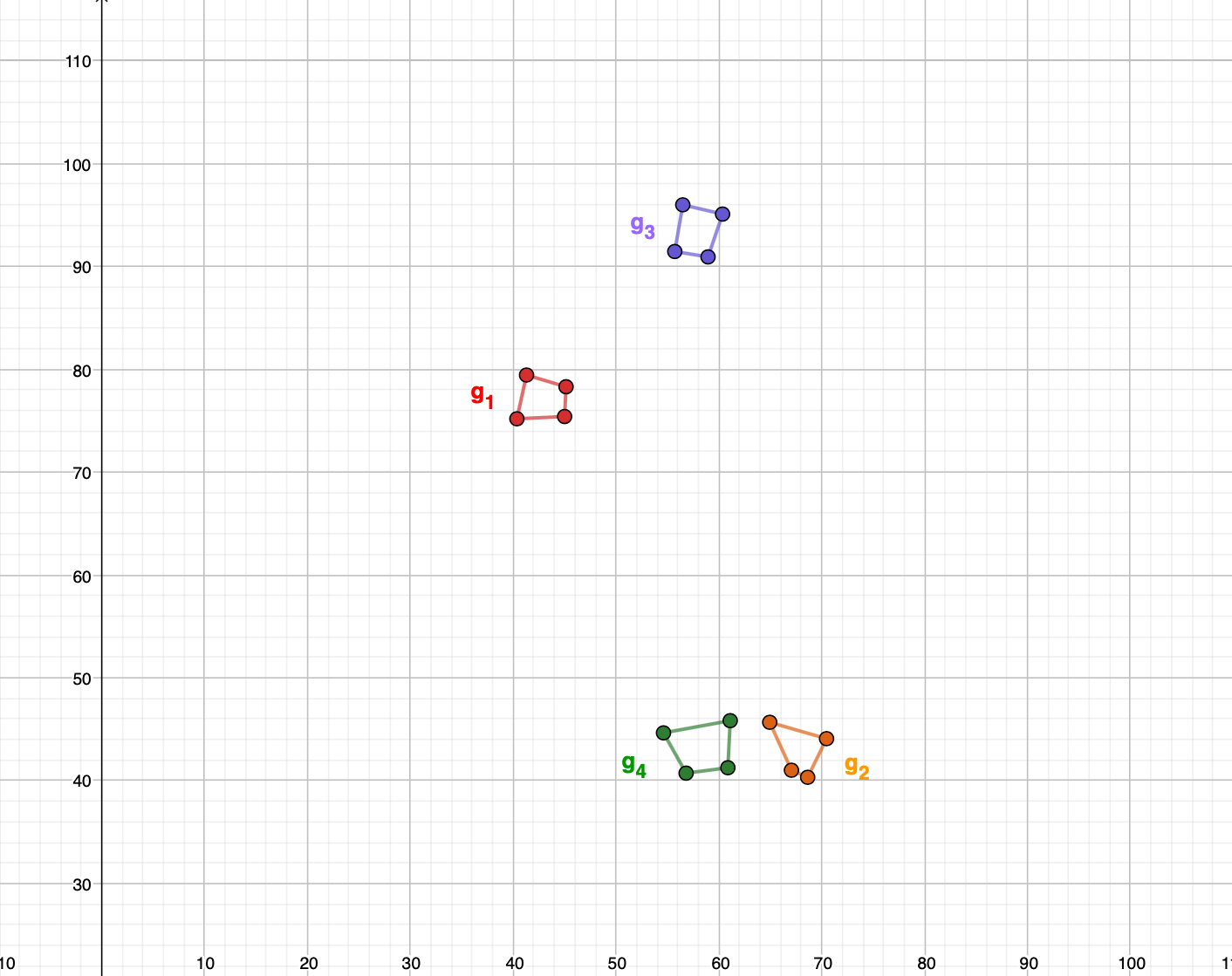}
\caption{Illustrative example \label{fig:example1}}
\end{figure}

\noindent
Figure \ref{fig:example1} shows an illustrative example consisting of four target planar graphs ($g_1$, $g_2$, $g_3$ and $g_4$) to be visited. We assume that their visits must be performed by a fleet of two drones supported by a mothership which starts from the origin $(0,0)$ and ends on the destination point $(100,0)$.

\begin{figure}[h!]
    \centering
    \subfloat[\centering a]{{\includegraphics[width=5cm]{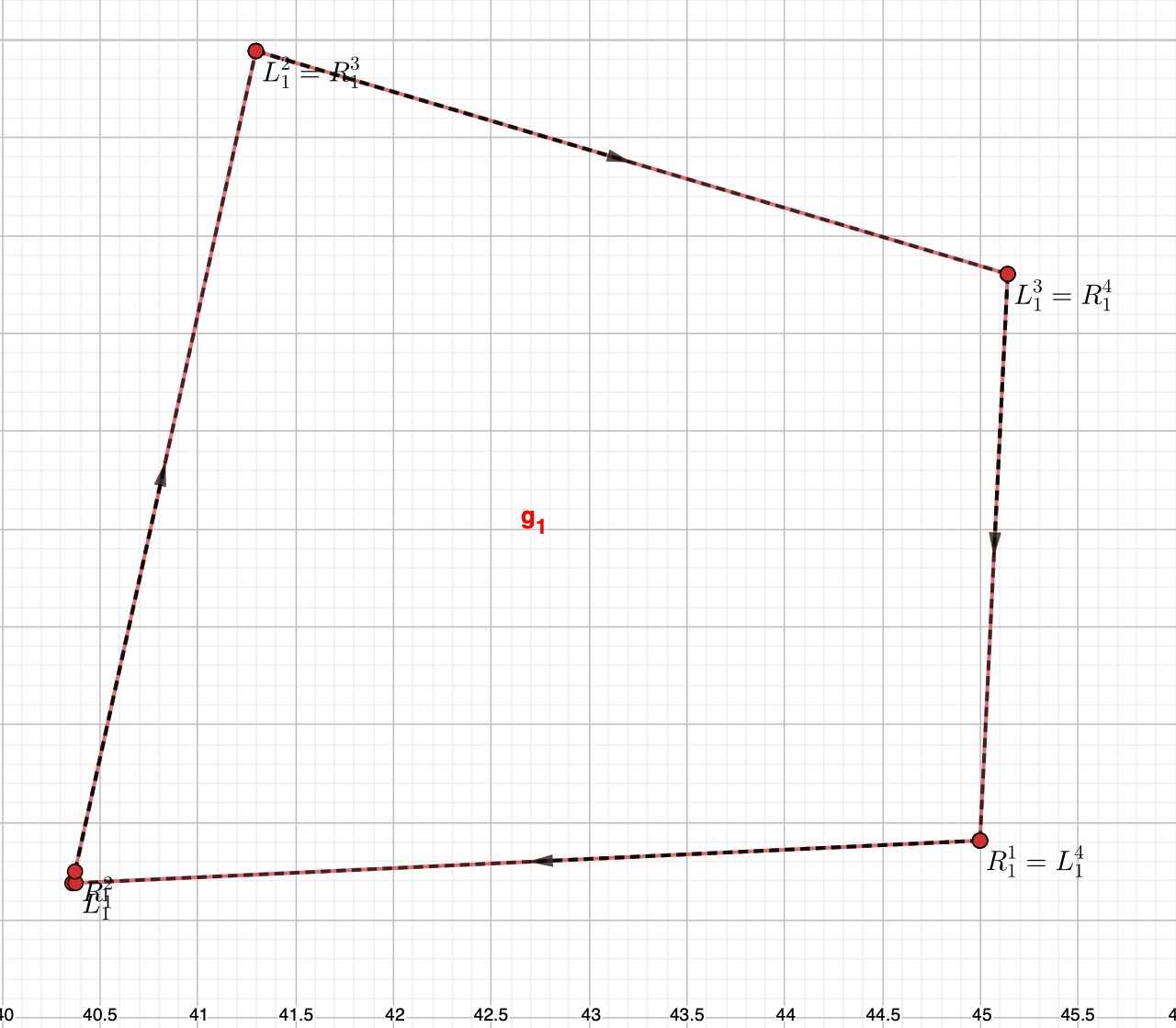} }}%
    \qquad
    \subfloat[\centering b]{{\includegraphics[width=5cm]{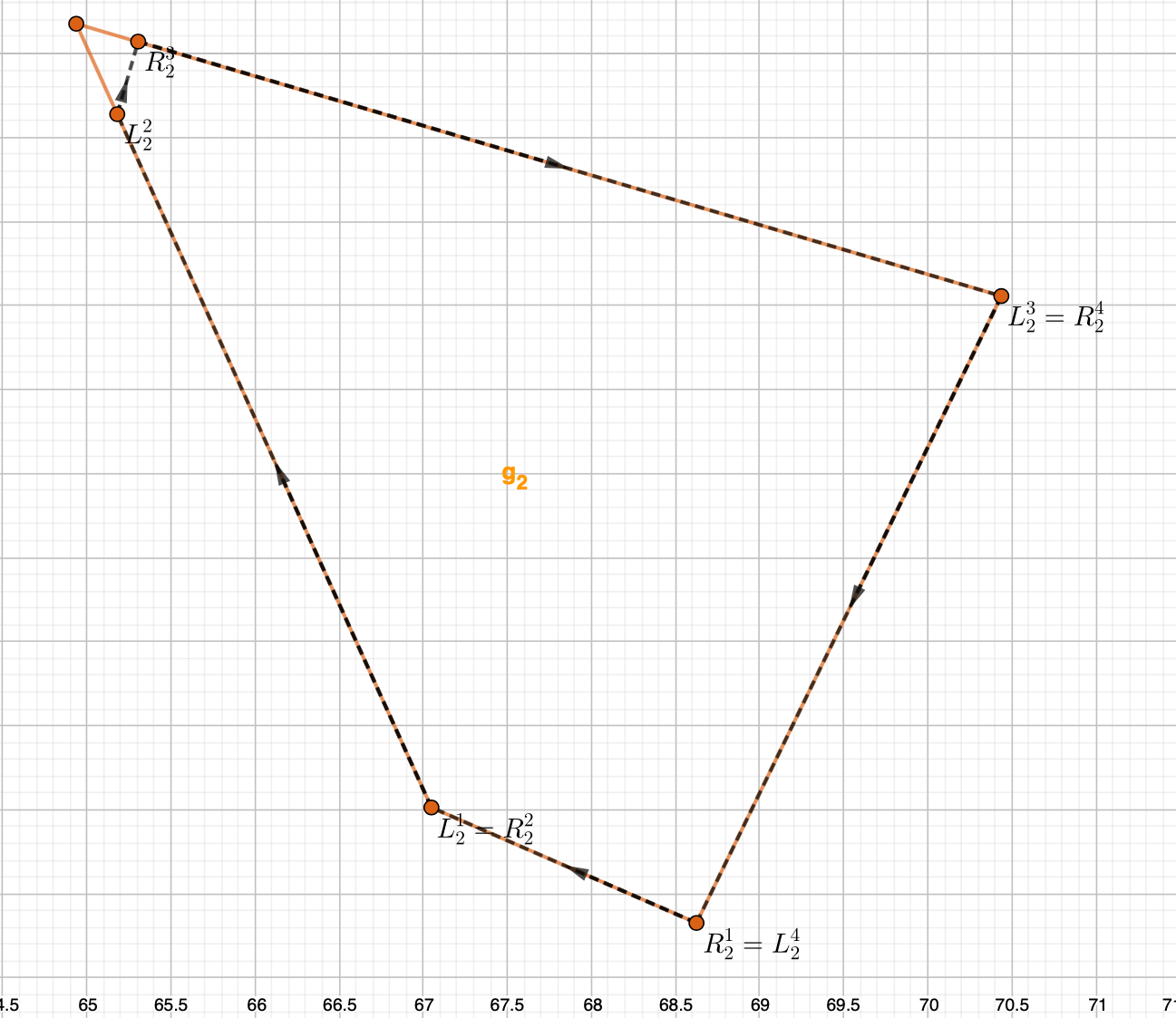} }}%
     \qquad
    \subfloat[\centering c]{{\includegraphics[width=5cm]{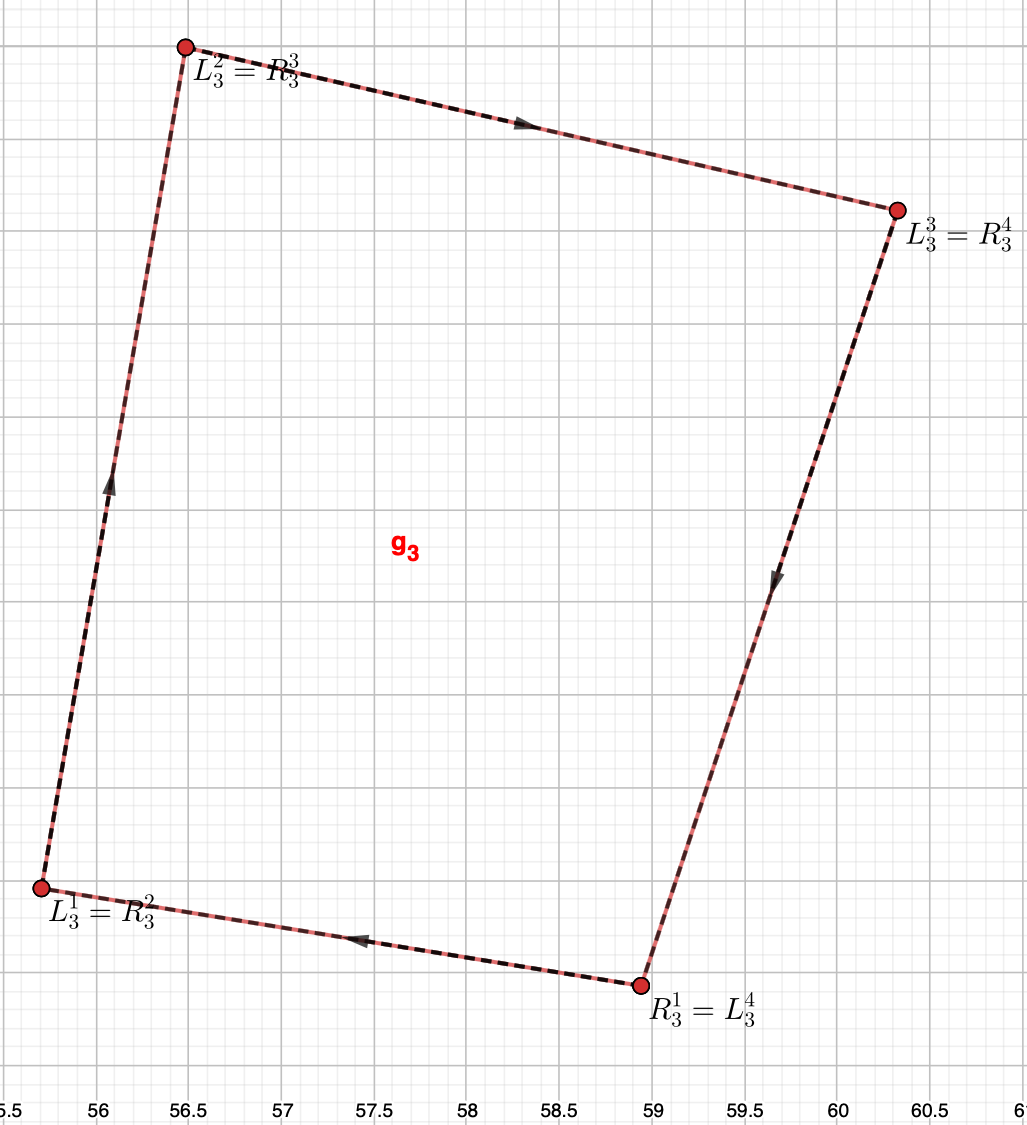} }}%
    \qquad
    \subfloat[\centering d]{{\includegraphics[width=5cm]{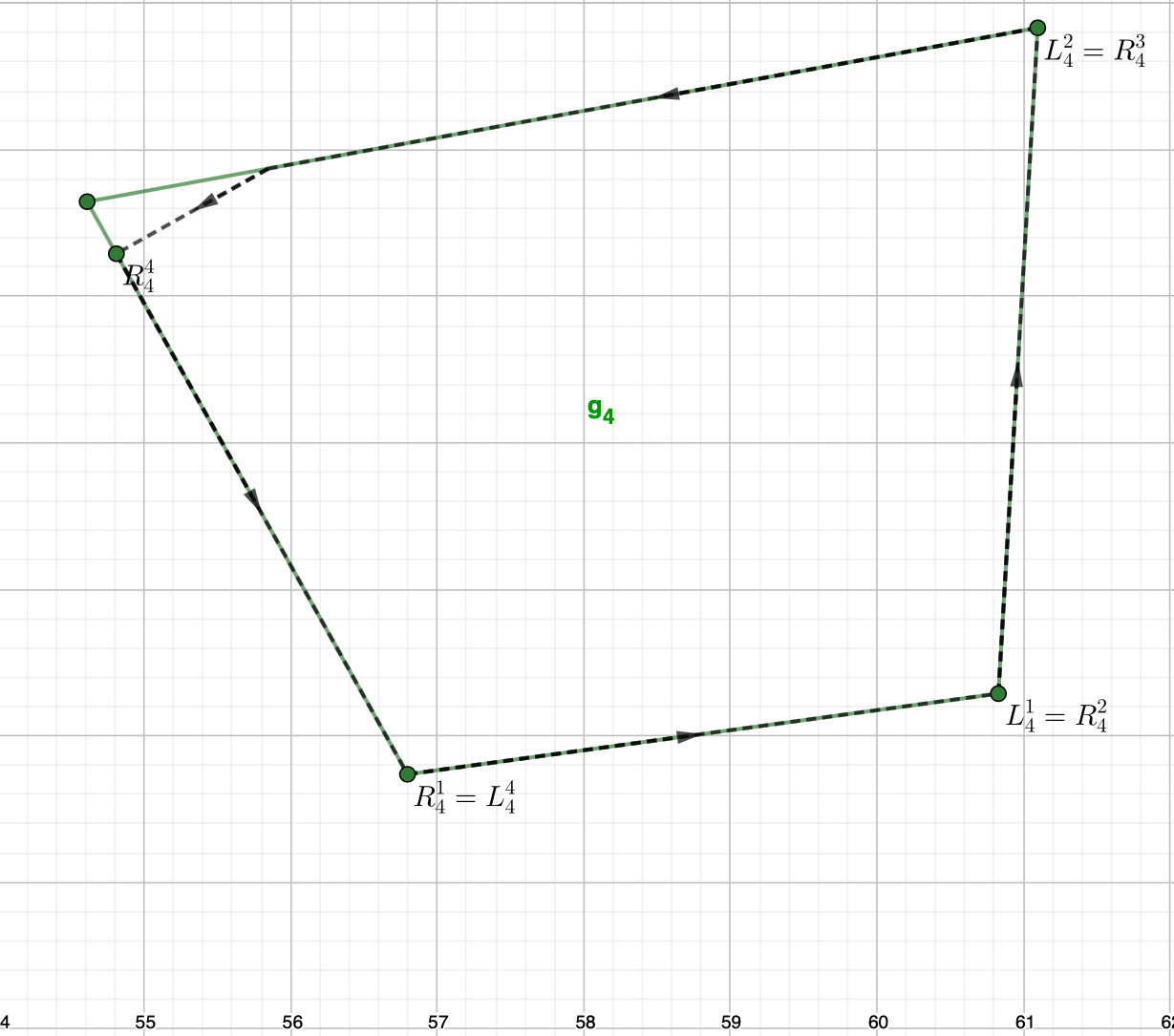} }}%
    \caption{STEP 1 for the illustrative example}%
    \label{fig:example2}%
\end{figure}

\noindent
Figure \ref{fig:example2} reports a zoom on each single target graph, showing the tour generated by STEP 1 of the heuristic procedure. A pair of points representing retrieving and launching points,  together with an arrow pointing the direction followed by the drone according with the order in which the edges are visited, are depicted on each edge.

\begin{figure}[h!]
    \centering
    \subfloat[\centering a]{{\includegraphics[width=5cm]{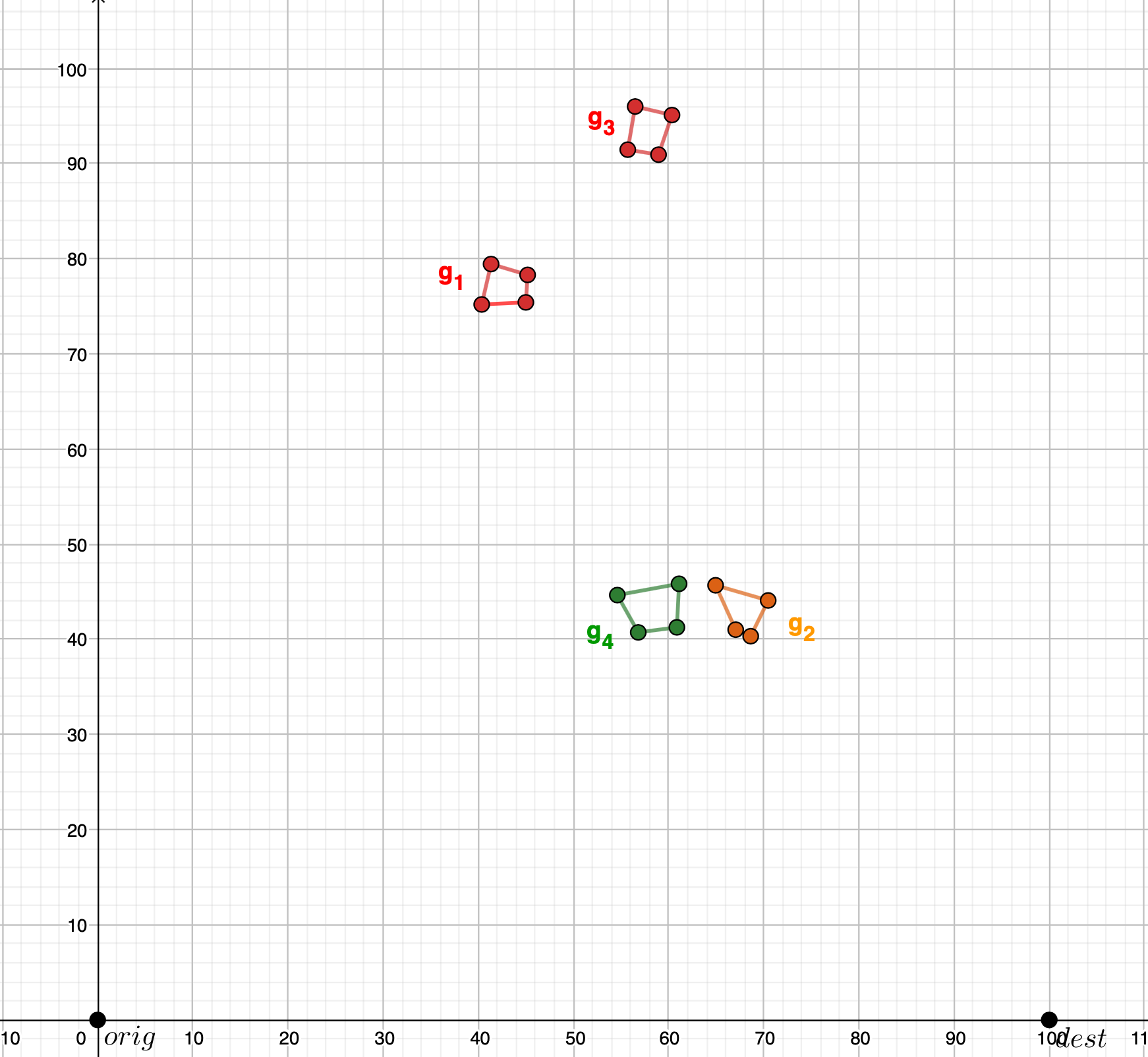} }}%
    \qquad
    \subfloat[\centering b]{{\includegraphics[width=5cm]{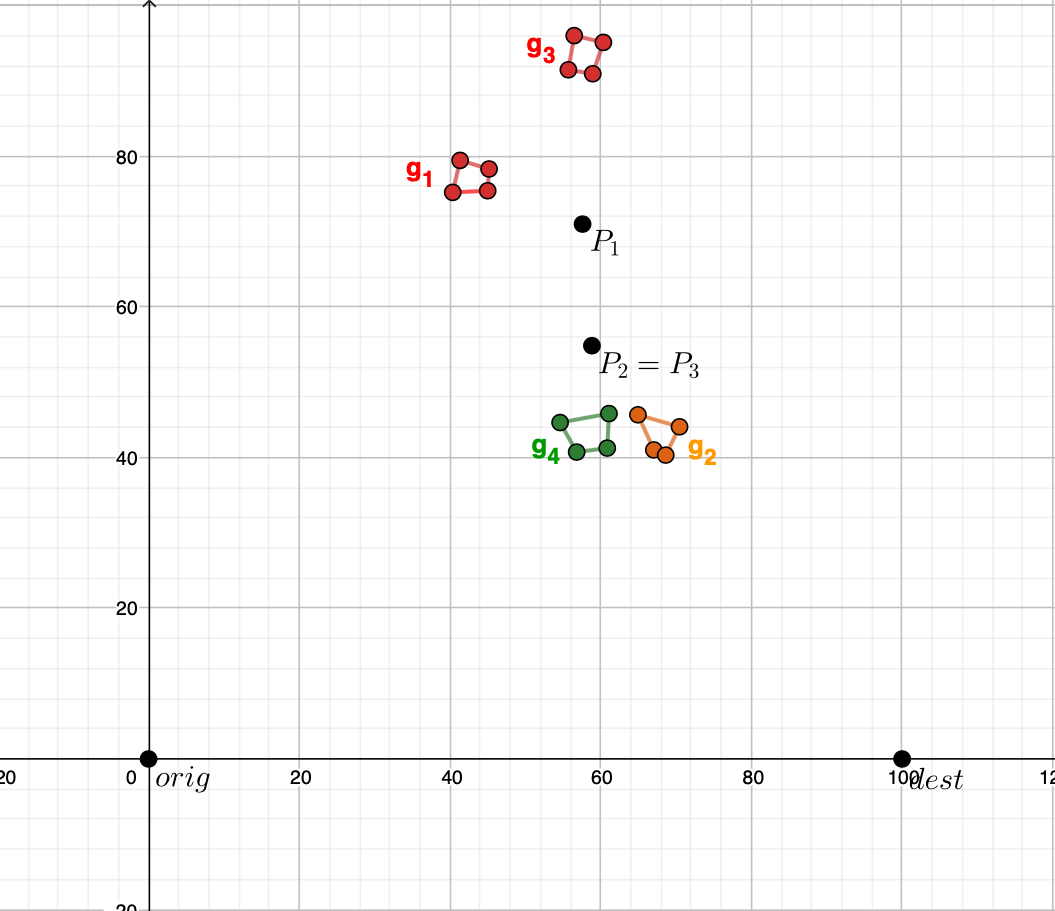} }}%
        \caption{[a] STEP 2, [b] STEP 3 for the illustrative example}%
    \label{fig:example3}%
\end{figure}

\noindent
By applying STEP 2 to this illustrative example, we obtain three clusters, as shown in Figure \ref{fig:example3}[a]. One cluster contains graphs $g_1$ and $g_3$ (in red), while graph $g_2$ and $g_4$ represent distinct clusters. The computation of the reference points of these clusters, according with STEP 3, produces the points $P_1$ and $P_2 = P_3$, as shown in Figure  \ref{fig:example3}[b]. Note that in this case the reference points for the two clusters containing respectively $g_2$ and $g_4$ coincide, and thus we have only two distinct reference points for the three clusters.


\begin{figure}[h!]
    \centering
    \subfloat[\centering a]{{\includegraphics[width=5cm]{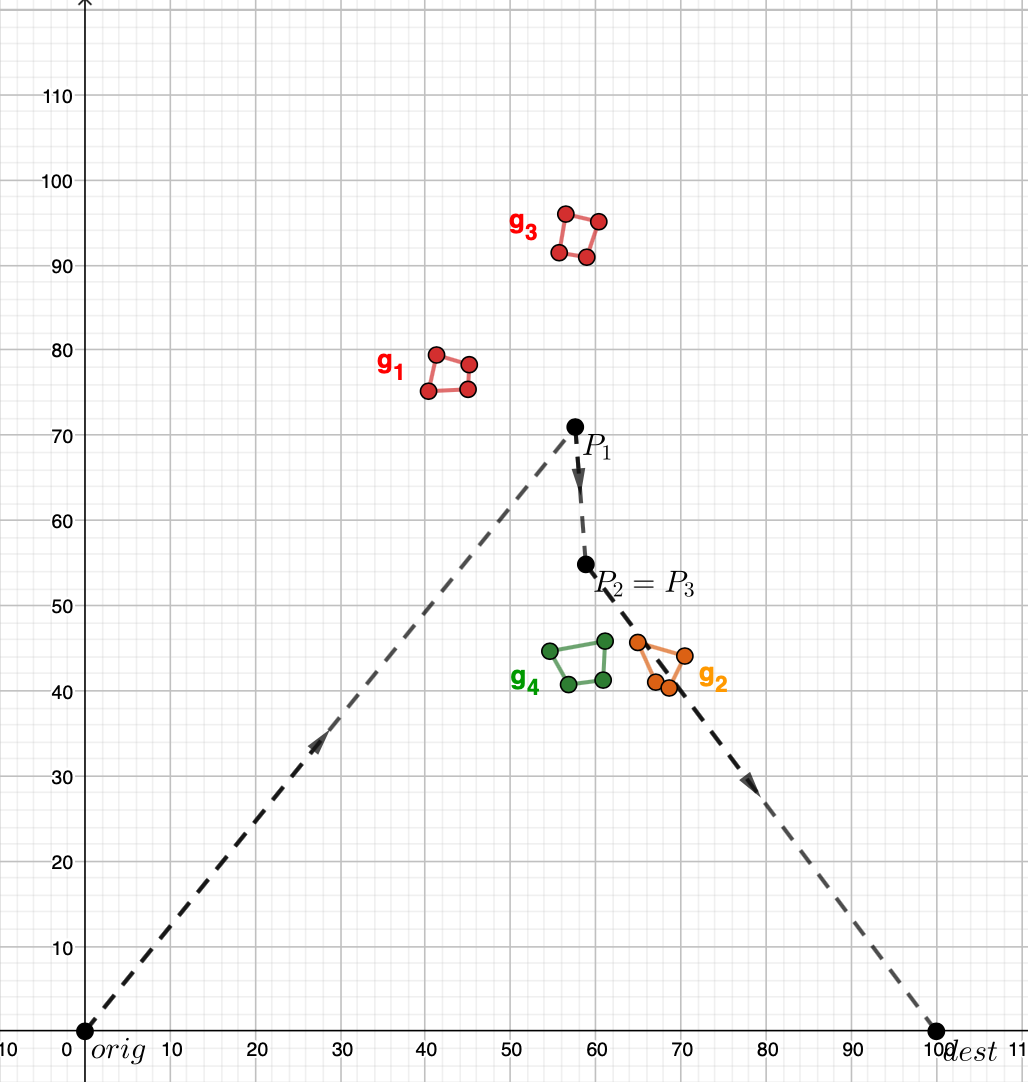} }}%
    \qquad
    \subfloat[\centering b]{{\includegraphics[width=5cm]{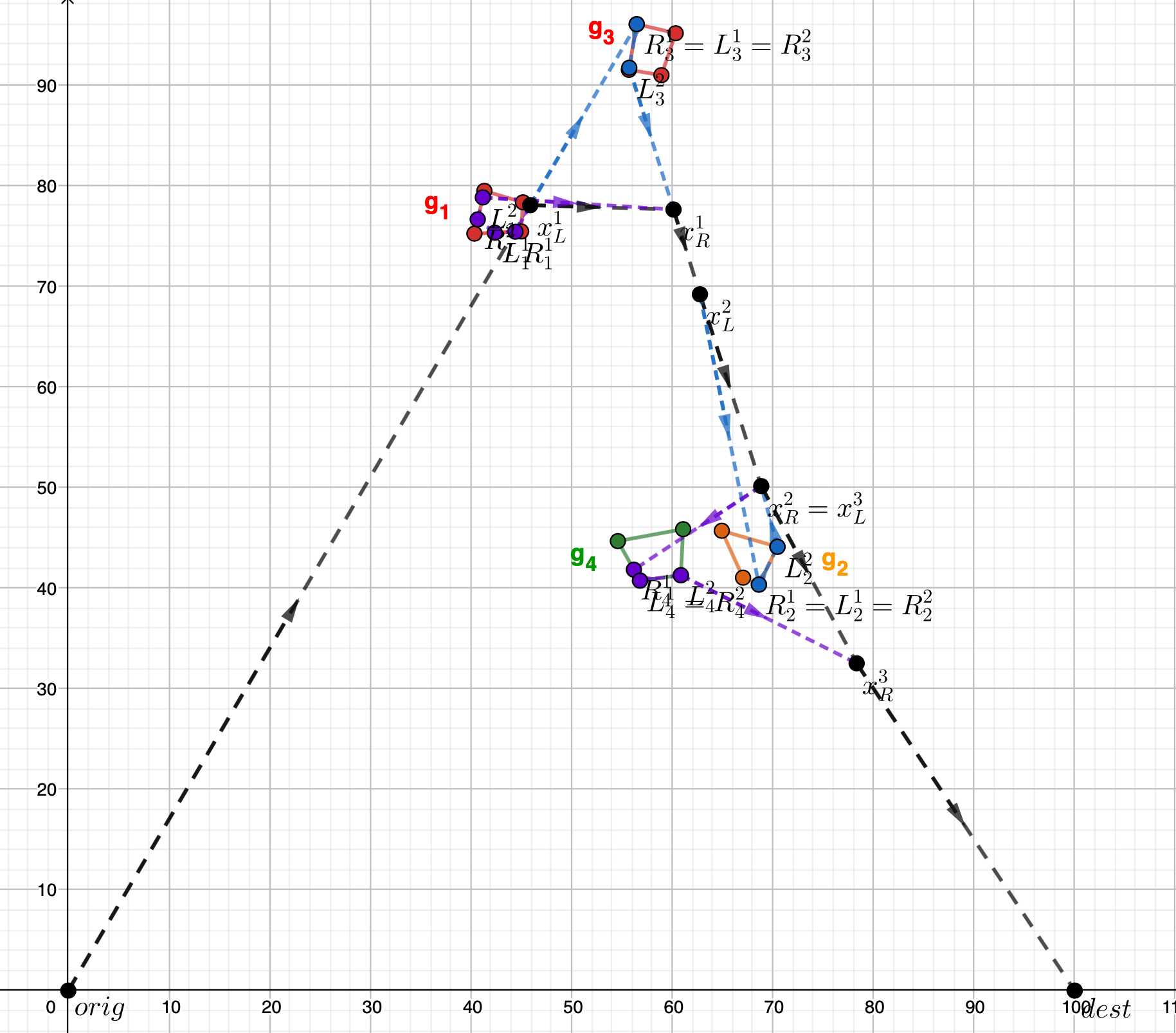} }}%
        \caption{[a] STEP 4, [b] STEP 5 for the illustrative example}%
    \label{fig:example4}%
\end{figure}

\noindent
STEP 4 of the heuristic procedure generates the tour of the mothership along the origin point, $P_1$, $P_2$ and the destination point, as shown in Figure \ref{fig:example4}[a]. This tour returns also the order in which the clusters are visited (and thus, also the order of visit to the target graphs) and this permits to set the values of the variables $u^{e_{g}td}$ and $v^{e_{g}td}$ of the \AMD\space model.\\
\noindent 
By providing the initial partial solution obtained by the values of the variables $u^{e_{g}td}$ and $v^{e_{g}td}$, STEP 5 solves the \AMD\space model and returns the final feasible solution shown in Figure \ref{fig:example4}[b]. From it we can observe that the sequence of visit of the target graphs does not change with respect to the one provided by STEP 4. The fleet of two drones first visits graphs $g_1$ and $g_3$ starting from the launching point $x^1_L$. Then, both drones are retrieved by the mothership at point $x^1_R$. The mothership moves to the point $x^2_L$  where one drone is launched for visiting graph $g_2$. Then the mothership reaches point $x^2_R$ to retrieve the drone and from the same point it launches the other drone for visiting graph $g_4$. Then, this drone is retrieved by the mothership at point $x^3_R$ before moving to the final destination point.

\begin{figure}[h!]
    \centering
    \subfloat[\centering a]{{\includegraphics[width=5cm]{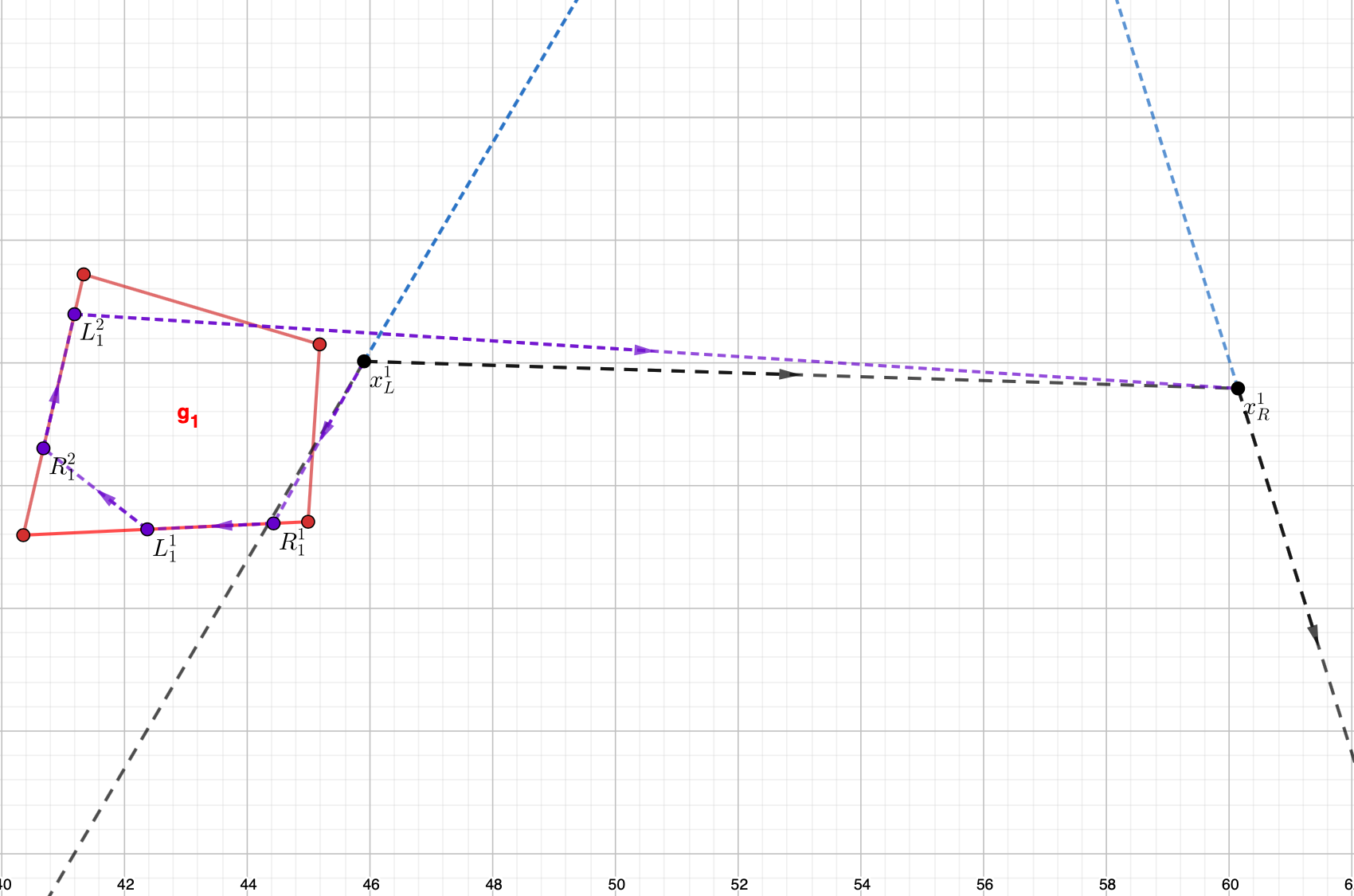} }}%
    \qquad
    \subfloat[\centering b]{{\includegraphics[width=5cm]{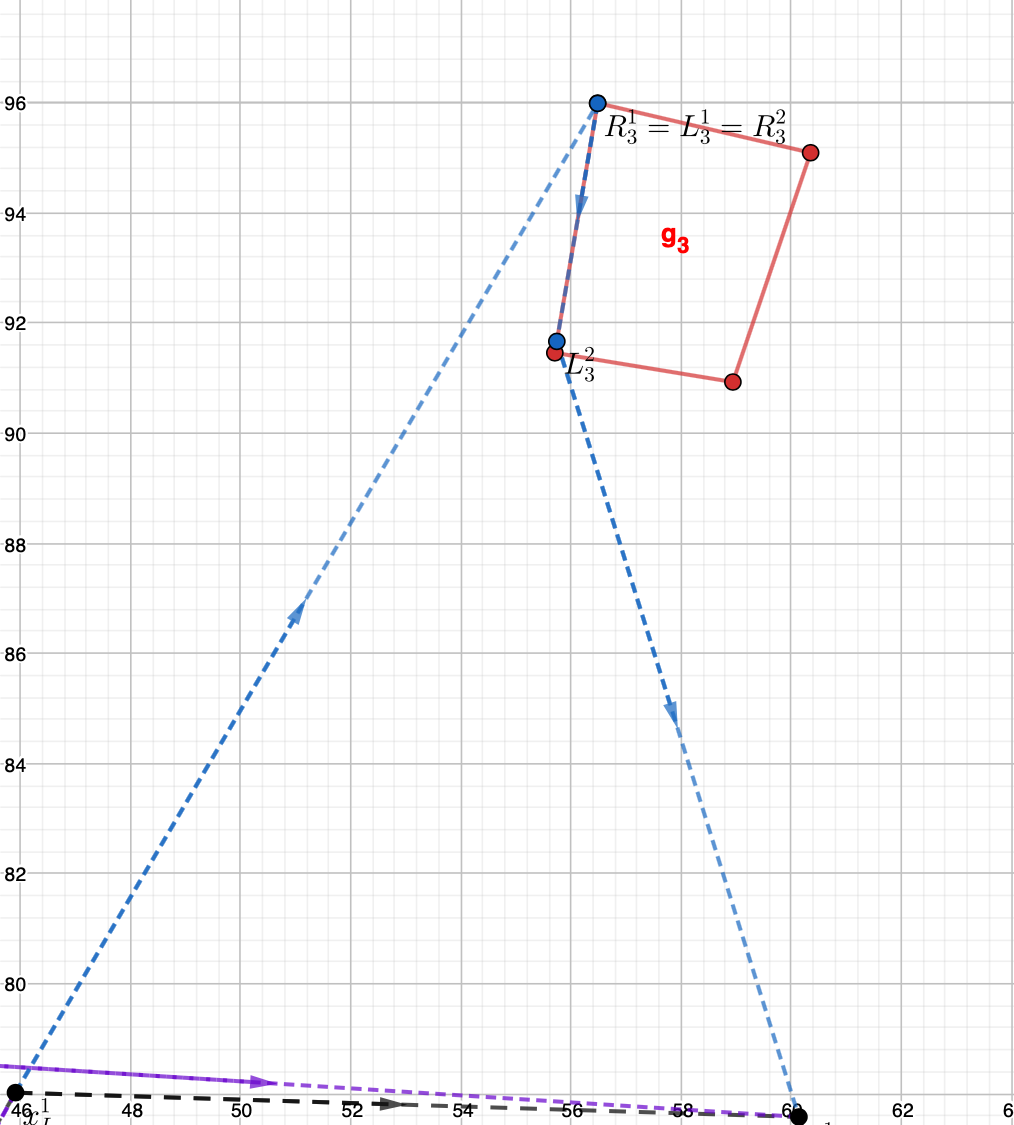} }}%
     \qquad
    \subfloat[\centering c]{{\includegraphics[width=5cm]{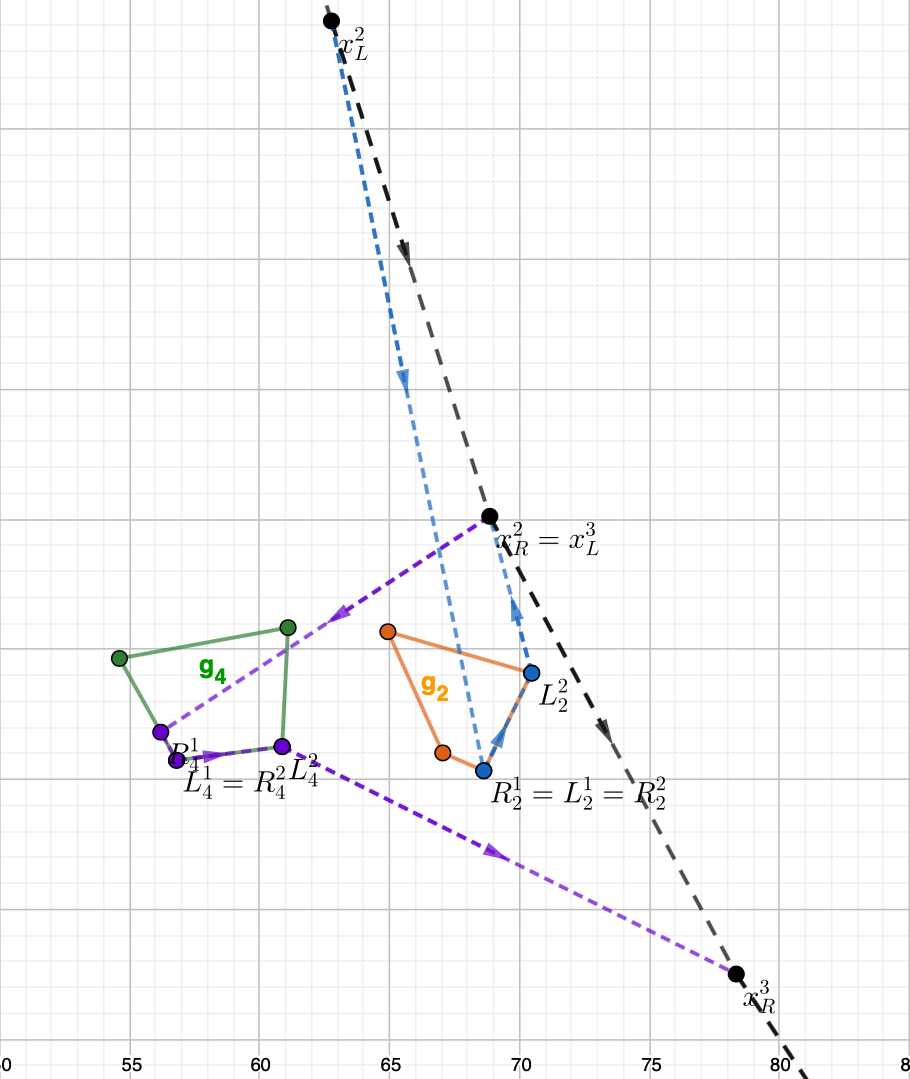} }}%
    \caption{Zoom on the tour on each target graph provided by STEP 5}%
    \label{fig:example5}%
\end{figure}

\noindent 
Focusing on each single target graph, Figure \ref{fig:example5} shows the zoom on the tours followed by the drones. For example, Figure \ref{fig:example5}[a] reports the one performed by the drone that visits graph $g_1$. The drone starts from the mothership at point $x^1_L$ and first visits the segment $\overline{R^1_1L^1_1}$. From point $L^1_1$ the drone moves to the second visited edge of graph $g_1$ visiting the segment $\overline{R^2_1L^2_1}$.
Finally, the drone leaves the graph from point $L^2_1$ and it is retrieved by the mothership at point $x^1_R$. Note that in this example the drones do not visit the full 100\% of each graph, but only a pre-specified percentage of each one of them.

\bigskip
\noindent
The reader may notice that the above algorithm can be also used to generate solutions for the model without synchronization  commented in Section \ref{amdasyn} since any solution of the synchronous model is also feasible for the asynchronous one. 

\newpage
\section{Experimental results}

\noindent
In this section we discuss the experimental results obtained testing the formulations presented in Section \ref{Form} and the matheuristic procedure proposed in Section \ref{Math} on a testbed of instances. In particular, we consider instances like the ones used in \cite{art:Amorosi2021} where the targets, to be visited by the drones, are represented by grid graphs.
This set consists of 5 instances of respectively 5 and 10 target graphs, with different cardinality of the set of nodes. More precisely, each instance is composed of 20$\%$ graphs with 4 nodes, 20$\%$ graphs of 6 nodes, 20$\%$ graphs of 8 nodes and 20$\%$ graphs of 10 nodes. Moreover, we assume that the velocity of the drones is twice that of the mothership and that a random percentage of each target graph, or of each of its edges, must be visited by the fleet of drones.\\
We consider in our experiments that the number of drones varies between 1 and 3 and that the drones
endurance (expressed as maximum time that the drone can operate when it is fully recharged) ranges between 20 and 60.
Table \ref{table:tab1} reports a summary of the characteristics of our instances.

\renewcommand{\arraystretch}{0.7}
\begin{table}[!h]
\caption{Instances parameter values}
\centering
\footnotesize
\begin{tabular}{l | c }
\hline
\# Targets & (5,10)\\
\hline
\# Drones &	(1,2,3)\\
\hline
\# Vertices & (4,6,8,10,12)\\
\hline
Drone endurance & (20,30,40,50,60)\\
\hline
$\%$ target (edge) & random variable
\end{tabular}
\label{table:tab1}
\end{table}

\noindent
Table \ref{table:tab2} reports the results obtained solving the \AMD model on the instances previously described, by adopting the commercial solver Gurobi. We consider the exact solution both providing and not providing an initial solution computed by the matheuristic described in Section \ref{Math}. More precisely, the first column of Table \ref{table:tab2} indicates the number of target graphs to be visited by the fleet of drones, the second column reports the capacity of the drones, the third column distinguishes between the visit of a percentage of each edge (e) and the percentage of each target graph (g). The fourth column reports the size of the fleet of drones. Thus, for each combination of the listed parameters, we summarize the average percentage gap of the solution obtained within the time limit set equal to 2 hours. We report respectively average percentage gap with initialization by the solution provided by the matheuristic, solution time, in seconds, of the matheuristic and average gap without initialization by the solution obtained by the matheuristic.\\
We can observe that the value of the average percentage gap ranges between a minimum of 66.9\% and a maximum of 97.43\%. This shows that the model is hard to be solved even on small size instances. Moreover, we can see that in most of the cases, the average percentage gap associated with the variant of the model consisting in visiting a given percentage of each edge, is higher than the one associated with the variant imposing to visit a given percentage of each target graph. Another thing that we can observe is that the average percentage gap increases with the number of drones and decreases with the drone endurance.\\
As regards the number of target graphs, we can see that by increasing it from 5 to 10, the exact method, without initialization by the solution obtained with the matheursitic, becomes even harder. Indeed, the red entries of the table mean that some instances could not find a feasible solution within the time limit (note that in the brackets we indicate the number of these instances). The number of not solvable instances increases with the number of drones. Moreover, for the minimum level of endurance, the exact solution of the model without initialization provided by the matheuristic, does not provide any solution, within the time limit, for instances with 10 graphs and 2 or 3 drones.\\
Considering the comparison with the exact method starting from the solution provided by the matheuristic, we can note that the values of average percentage gap are very close to the ones related to the exact solution method without initialization. Thus the initialization does not speed up the convergence of the solver. However, we can see that the matheuristic is always able to find a feasible solution of the problem, even for the cases in which the solver is not.\\
Moreover, the average solution times of the matheuristic range between a minimum of 37 seconds to a maximum of 3 minutes. They increase with the drone capacity for the variant of the model in which a given percentage of each edge must be visited, while they decrease by increasing the number of drones for the variant of the model in which a given percentage of each target graph must visited. By increasing the number of target graphs from 5 to 10, the average solution times of the matheuristics become more than double for both model variants.
Summing up, the results obtained show that the exact solution method given by solving the formulation is very challenging even for small size instances. However, exploiting it, the matheuristic is able to provide  solutions for all instances rather quickly.

\begin{table}[h!]
\centering
\caption{Comparison between exact solution with and without initialization by the matheuristic solution}
\label{table:tab2}
\resizebox{\textwidth}{!}{%
\begin{tabular}{|c|c|c|ccc|ccc|ccc|}
\hline
 &  &  & \multicolumn{9}{c|}{\textbf{\# drones}} \\ \cline{4-12} 
 &  &  & \multicolumn{3}{c|}{1} & \multicolumn{3}{c|}{2} & \multicolumn{3}{c|}{3} \\ \cline{4-12} 
\multirow{-3}{*}{$\bm{|\mathcal G|}$} & \multirow{-3}{*}{$\bm{N^d}$} & \multirow{-3}{*}{\textbf{v.t.}} & \% Gap (i) & TimeH & \% Gap (wi) & \% Gap (i) & TimeH & \% Gap (wi) & \% Gap (i) & TimeH & \% Gap (wi) \\ \hline
 &  & e & 82,63 & 61,56 & 81,70 & 91,57 & 63,80 & 90,61 & 93,06 & 60,87 & 90,93 \\ \cline{3-3}
 & \multirow{-2}{*}{20} & g & 79,09 & 44,97 & 79,63 & {\color[HTML]{333333} 89,03} & 37,32 & 91,85 & 94,00 & 39,05 & 95,80 \\ \cline{2-12} 
 &  & e & 82,70 & 65,21 & 80,17 & 85,14 & 64,41 & 82,21 & 91,9 & 63,34 & 90,12 \\ \cline{3-3}
 & \multirow{-2}{*}{30} & g & 75,80 & 55,77 & 71,19 & 84,36 & 44,36 & 88,27 & 91,02 & 44,59 & 91,39 \\ \cline{2-12} 
 &  & e & 80,94 & 68,81 & 77,98 & 83,44 & 64,80 & 82,16 & 91,24 & 63,19 & 86,25 \\ \cline{3-3}
 & \multirow{-2}{*}{40} & g & 74,47 & 43,92 & 73,46 & 81,21 & 38,27 & 84,35 & 85,34 & 37,51 & 89,63 \\ \cline{2-12} 
 &  & e & 76,87 & 66,67 & 74,41 & 81,12 & 63,86 & 79,57 & 85,11 & 63,51 & 86,16 \\ \cline{3-3}
 & \multirow{-2}{*}{50} & g & 70,58 & 43,42 & 66,90 & 80,96 & 43,98 & 88,84 & 80,49 & 44,35 & 82,81 \\ \cline{2-12} 
 &  & e & 76,39 & 67,78 & 71,61 & 81,63 & 66,08 & 79,84 & 83,82 & 64,40 & 82,06 \\ \cline{3-3}
\multirow{-10}{*}{5} & \multirow{-2}{*}{60} & g & 78,17 & 44,69 & 72,79 & 79,35 & 40,63 & 86,55 & 81,74 & 50,01 & 84,66 \\ \hline
 &  & e & 82,56 & 137,93 & 84,91 & 92,30 & 128,53 & - & 94,73 & 124,44 & - \\ \cline{3-3}
 & \multirow{-2}{*}{20} & g & 81,00 & 119,20 & {\color[HTML]{FE0000} 84,08 (2)} & 89,88 & 83,50 & {\color[HTML]{FE0000} 96,64 (2)} & 96,44 & 70,00 & {\color[HTML]{FE0000} 97,43 (3)} \\ \cline{2-12} 
 &  & e & 80,60 & 159,00 & 80,93 & 87,11 & 132,15 & {\color[HTML]{FE0000} 87,58 (3)} & 94,56 & 127,35 & {\color[HTML]{FE0000} 92,85 (2)} \\ \cline{3-3}
 & \multirow{-2}{*}{30} & g & 79,93 & 132,67 & {\color[HTML]{FE0000} 82,70 (1)} & 86,32 & 80,29 & {\color[HTML]{FE0000} 86,13 (3)} & 91,12 & 76,72 & {\color[HTML]{FE0000} 89,74 (1)} \\ \cline{2-12} 
 &  & e & 79,05 & 191,37 & 78,07 & 85,11 & 131,26 & 84,33 & 91,88 & 132,10 & {\color[HTML]{FE0000} 88,61 (1)} \\ \cline{3-3}
 & \multirow{-2}{*}{40} & g & 80,23 & 115,00 & 79,64 & 87,31 & 68,39 & {\color[HTML]{FE0000} 84,57 (3)} & 96,09 & 69,40 & {\color[HTML]{FE0000} 91,86 (1)} \\ \cline{2-12} 
 &  & e & 81,49 & 188,32 & 77,81 & 87,72 & 134,01 & {\color[HTML]{FE0000} 85,51 (1)} & 92,68 & 132,82 & {\color[HTML]{FE0000} 90,79 (3)} \\ \cline{3-3}
 & \multirow{-2}{*}{50} & g & 79,92 & 87,23 & 80,38 & 82,80 & 66,14 & {\color[HTML]{FE0000} 84,00 (3)} & 92,48 & 64,94 & {\color[HTML]{FE0000} 91,96 (2)} \\ \cline{2-12} 
 &  & e & 83,79 & 155,27 & 81,57 & 85,91 & 131,94 & {\color[HTML]{FE0000} 82,96 (2)} & 92,24 & 130,11 & {\color[HTML]{FE0000} 86,58 (3)} \\ \cline{3-3}
\multirow{-10}{*}{10} & \multirow{-2}{*}{60} & g & 77,57 & 97,89 & 78,46 & 86,94 & 76,53 & {\color[HTML]{FE0000} 88,29 (2)} & 94,31 & 69,53 & {\color[HTML]{FE0000} 92,23 (3)} \\ \hline
\end{tabular}%
}
\end{table}

\begin{figure}[h!]
\includegraphics[width=\linewidth]{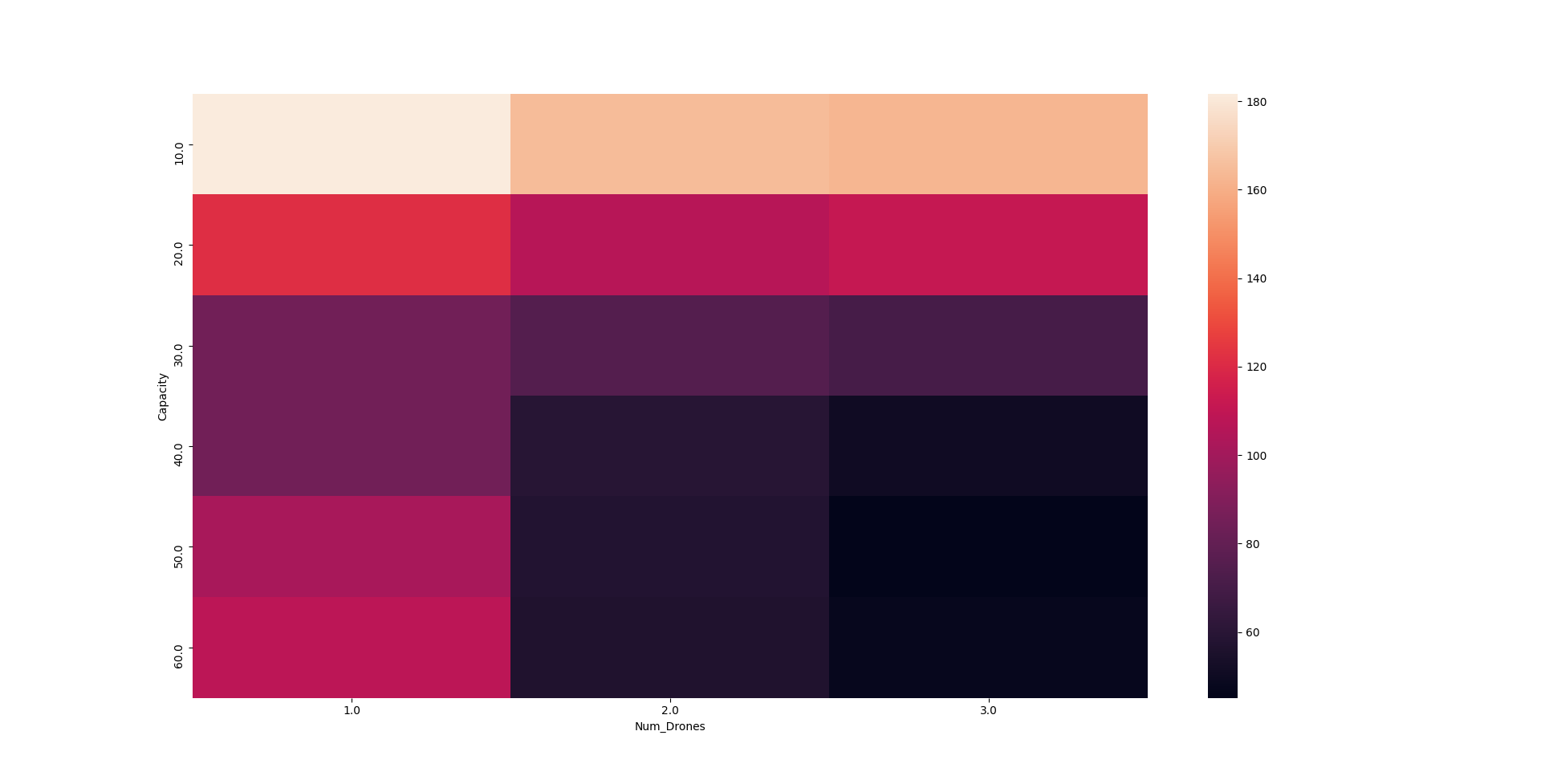}
\caption{Heatmap of objective function values depending on number of drones and drone capacities. The darker the color intensity the smaller the objective value. \label{fig:heatmap}}
\end{figure}

\noindent
We show in Figure \ref{fig:heatmap} the relationship between the objective function values of the problem and the number of available drones and their capacity. This figure reports the average objective values of all the instances with three target graphs varying the number of drones in $\{1,2,3\}$ and their endurance (capacity) in $\{10,20,30,40,50,60\}$. The darker the color intensity the smaller the objective value. As expected, our experiment confirms that both, a greater number of drones and larger endurance reduce total length of the mothership route.

\noindent

\noindent

\noindent
\noindent

\noindent

\noindent

\noindent

\section{Case Study}

\begin{figure}[h!]
\centering
\includegraphics[width=0.6\linewidth]{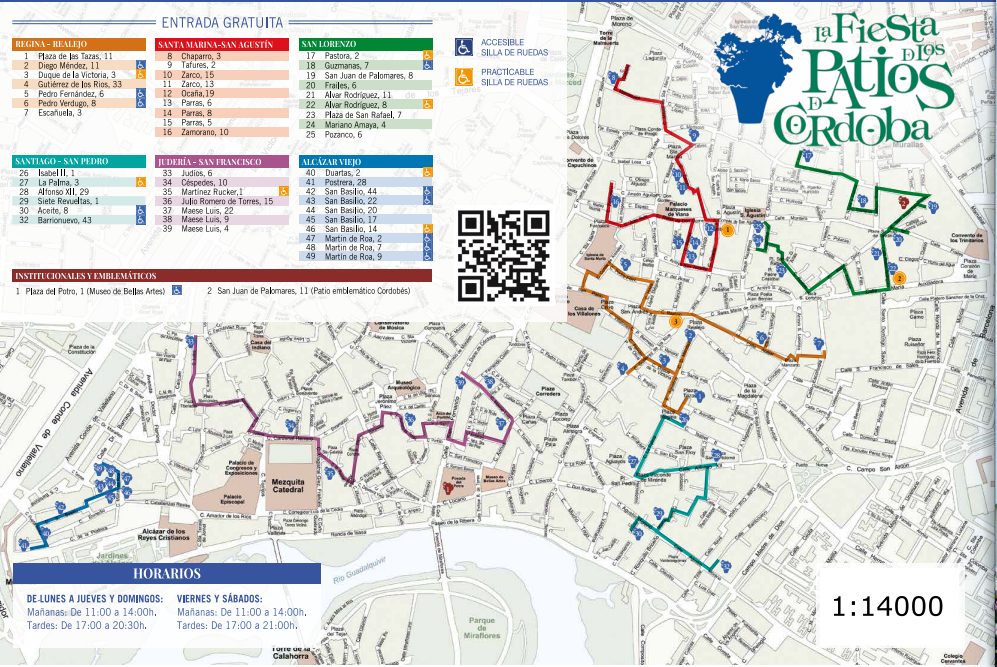}
\caption{Map of the Courtyards Festival in Cordoba. \label{fig:mapPF}}
\end{figure}

\noindent
In this section we describe a realistic application of the system studied in this paper to perform surveillance operations. Considering the current COVID-19 restrictions, we focus on the problem of preventing and identifying possible concentrations of people during events such as popular or religious festivals. In particular we consider the Courtyards Festival of Cordoba (\url{https://patios.cordoba.es/es/}). This is a social event that takes place every year in the city of Cordoba, Spain, during the first two weeks of May. Courtyard’s owners decorate their houses with many flowers trying to win the award that is offered by the Town Hall. During this competition a festival runs in parallel with a number of artistic performances along six different paths located in different areas in the city as shown in Figure \ref{fig:mapPF}.
In the pandemic context, to monitor the situation to avoid concentration of people, we propose to apply a system consisting of one helicopter and a fleet of three drones.
This kind of system has been proved successfully and has been already applied in the military field by the US Army in order to leave the helicopter to the edge of dangerous airspace and release drones, which will then penetrate enemy territory and send back intelligence, surveillance and reconnaissance information (see \cite{FG}).
In our application the reason to adopt a similar system is the possibility to inspect simultaneously and in real time different paths also reducing the risk of flying the helicopter over populated areas and the cost for moving the helicopter by minimizing the total length of its tour.

\begin{figure}[h!]
\centering
\includegraphics[width=0.6\linewidth]{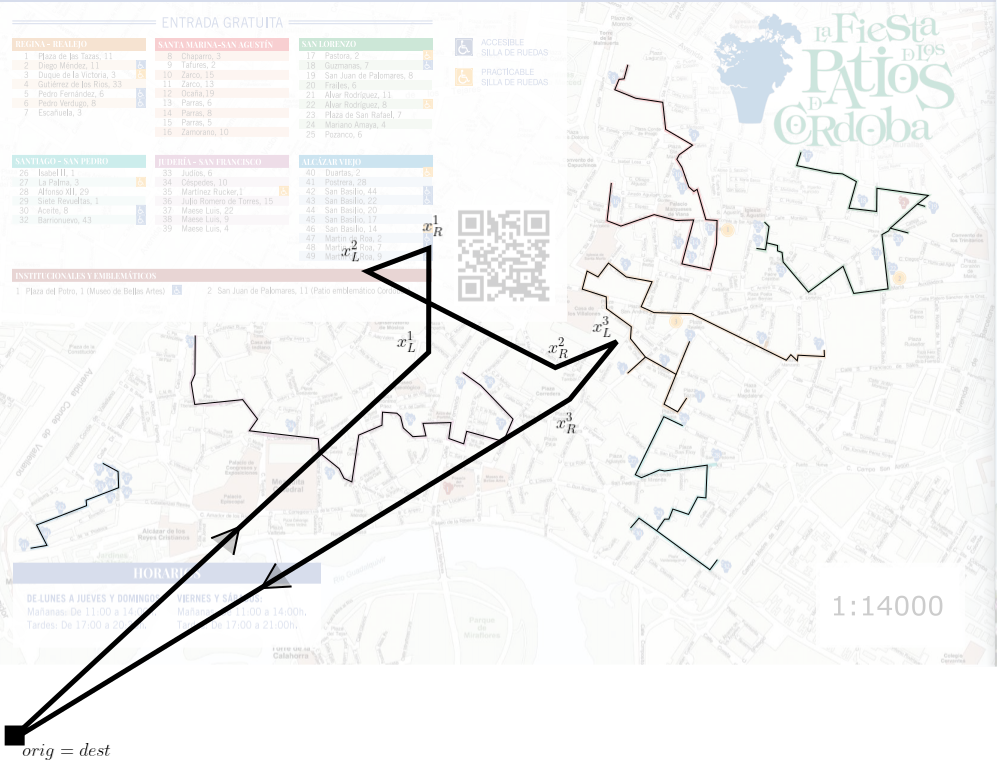}
\caption{Mothership tour. \label{fig:Mtour}}
\end{figure}

\begin{figure}[h!]
\centering
\includegraphics[width=0.6\linewidth]{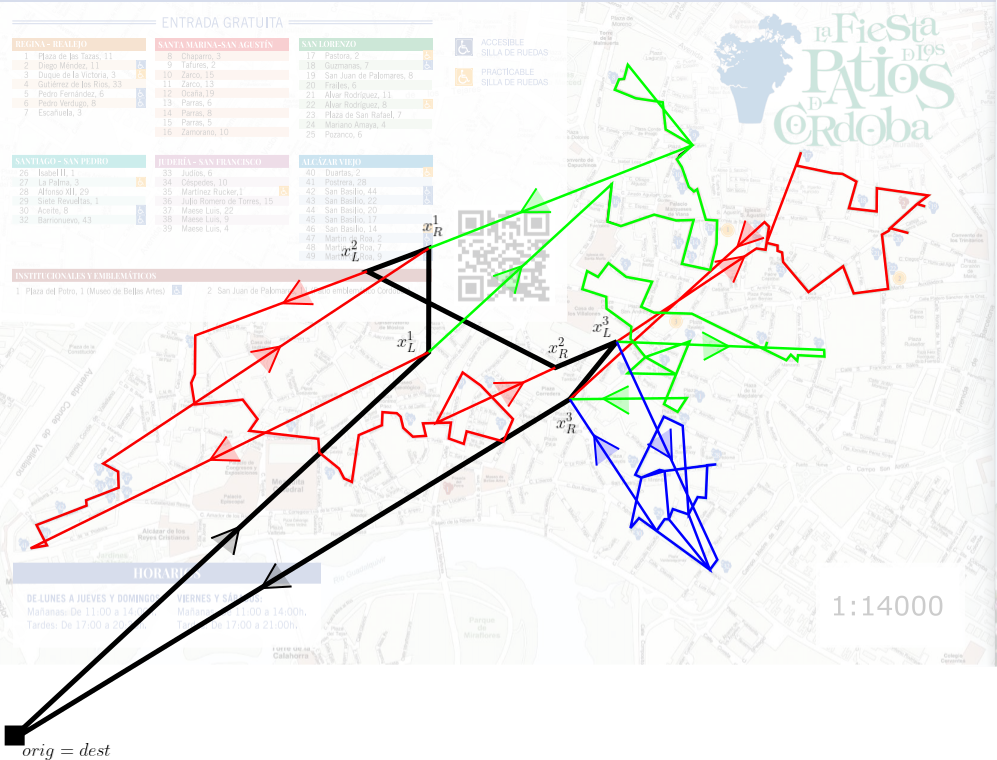}
\caption{The complete solution. \label{fig:tourD}}
\end{figure}

\noindent
We run the model presented in Section \ref{Form} on this scenario starting from the initial solution provided by the matheuristic, where the 6 coloured paths reported in the map of Figure \ref{fig:mapPF} represent the 6 target graphs to be visited, in this case inspected, by the fleet of drones. In addition, we suppose that the drones' speed is 100 km/h while that of the helicopter is  50 km/h aiming to minimize costs.
Moreover, we assume that the fleet is composed by three drones with an endurance equal to 2 hours, and we impose that each target graph must be fully visited (inspected).  As we can see from Figure \ref{fig:Mtour}, the origin of the mothership tour coincides with the destination and it is located in an area of the city where it is possible to assume the take-off and landing of an helicopter. Figure \ref{fig:Mtour} reports the tour followed by the helicopter in the solution obtained within the time limit of 2 hours sets to solve the model and with a percentage gap equal to 83\%. We can observe that the helicopter, starting from the origin, flies to the point $x_L^1$ that is the first launching point and then flies along the edge connecting $x_L^1$ with $x_R^1$, that is the first rendezvous point. Next, the helicopter flies to $x_L^2$ for launching the second drones' mission that are retrieved at point $x_R^2$. The third and last mission starts from $x_L^3$ and ends at point $x_R^3$ from where the helicopter goes back to the final destination.\\
Figure \ref{fig:tourD} shows the tour followed by the three drones for inspecting the six paths. In particular, one drone, in red, starts from $x_L^1$ for visiting the path of "Alcazar Viejo". From the same point a second drone, in green, starts for visiting the path of "Santa Maria-San Agustin". Both drones end their first mission at point $x_R^1$, where they are retrieved by the helicopter. Then, the helicopter flies to point $x_L^2$ where only one drone, the red one, starts its second mission to visit the path of "Juderia-San Francisco". In the meanwhile the helicopter, containing the other two drones, flies to point $x_R^2$ where it retrieves the previously launched drone. The last mission involves all the three drones that are launched from point $x_L^3$. The first drone, the red one, visits the path of "San Lorenzo", the second, the green one, visits the path of "Regina-Realejo" and the third, the blue one, visits the path of "Santiago-San Pedro".
All the three drones are retrieved by the helicopter at point $x_R^3$ and after that the helicopter goes back to the destination. 
The total distance travelled by the helicopter is equal to 11.27 km.

All details of this case study, including maps coordinates, .lp models and solutions can be found in \cite{Puerto2021}.

\section{Concluding remarks}
\noindent
This papers has analyzed the coordination problem that arises between a mothership vehicle and a fleet of drones that must coordinate their routes to minimize the total distance travelled by the mothership while visiting a set of targets modeled by graphs. We have presented exact formulations of the problem, for its synchronized and not synchronized versions. They are mixed integer non-linear programming models. Moreover, we presented valid inequalities for them. \\
Our computational results show that the considered problem is very challenging to solve even on small to medium size instances. For that reason, additionally, we have proposed a matheuristic algorithm that provides acceptable feasible solutions in very short computing time; so that it is a good alternative to the exact method. We report extensive computational experiments on randomly generated instances. Moreover, we present a case study related to inspection activities in the context of COVID-19 restrictions. We show the application of the system described in this paper in the framework of the Courtyard Festival in the city of Cordoba, by illustrating the solution obtained by adopting the problem formulation, in its synchronized version, and its solution by means of the initialization provided by the proposed matheuristic.\\
\noindent
The formulation and algorithms proposed in this paper can be seen a first building block to handle coordination of systems given by a base vehicle and drones. Further research in this topic must focus on finding faster and more accurate algorithms able to solve larger instances. Moreover, it is also challenging to model more complex operations  allowing that drones can visit more than one target per trip. These problems being very interesting are beyond the scope of the present paper and will be the focus of a follow up research line.
\section*{Acknowledgements}
This research has been partially supported by Spanish Ministry of Education and Science/FEDER grant number  MTM2016-74983-C02-(01-02), and projects Junta de Andalucia P18-FR-1422, FEDER-US-1256951, CEI-3-FQM331 and  \textit{NetmeetData}: Ayudas Fundaci\'on BBVA a equipos de investigaci\'on cient\'ifica 2019.


\bibliographystyle{elsarticle-num-names}

\bibliography{bibliography2.bib}

\end{document}